\documentclass{article}

\usepackage{amsmath, latexsym, amsfonts, amssymb, amsthm, amscd}

\newtheorem {lemme} {Lemma} [section]
\newtheorem {theoreme} {Theorem} [section]
\newtheorem {proposition} {Proposition} [section]
\newtheorem {corollaire} {Corollary} [section]
\newtheorem {remarque} {Remark} [section]

\newcommand{\R}{\mathbb {R}}
\newcommand{\C}{\mathbb {C}}

\newcommand{\1}{1\!\!{\sf I}}

\numberwithin{equation}{section}

\newcommand{\vers}{\mathop{\longrightarrow }}
\title{Exact separation phenomenon for the eigenvalues of large Information-Plus-Noise type matrices\\Application to spiked models}

\author{M. Capitaine\thanks{CNRS, Institut de Math\'ematiques de Toulouse, 
Equipe de Statistique et Probabilit\'es,  F-31062 Toulouse Cedex 09. 
E-mail: mireille.capitaine@math.univ-toulouse.fr }}

\date{}
\begin{document}
\maketitle
\begin{abstract}
We consider large Information-Plus-Noise type matrices  of the form  $M_N=( \sigma \frac{ X_N}{\sqrt{N}}+A_N)(\sigma \frac{ X_N}{\sqrt{N}}+A_N)^*$ where $X_N$ is an $n \times N$ ($n\leq N)$ matrix consisting of independent standardized complex entries, $A_N$ is  an $n \times N$ nonrandom matrix
and $\sigma>0$. As $N$ tends to infinity, if $n/N \rightarrow c\in ]0,1]$ and if the empirical spectral measure of  $A_N A_N^*$ converges weakly  to some compactly supported probability distribution $\nu \neq \delta_0$, Dozier and Silverstein established 
in \cite{DozierSilver} that almost surely the empirical spectral measure  of $M_N$ converges weakly  towards a nonrandom distribution $\mu_{\sigma,\nu,c}$. In   \cite{BaiSilver}, Bai and Silverstein proved, under certain assumptions on the model, that for some closed interval in $]0;+\infty[$
outside the support of $\mu_{\sigma,\nu,c}$ satisfying some conditions involving $A_N$,  almost surely, no eigenvalues of $M_N$ will appear in this interval for all $N$ large. In this paper,   we carry on with the study of the  support of the limiting spectral measure previously investigated in \cite{DozierSilver2} and later
in  \cite{VLM, LV} and we
show that, under almost  the same assumptions  as in  \cite{BaiSilver}, there is an exact separation phenomenon between the spectrum of $M_N$ and the spectrum of $A_NA_N^*$: 
 to a gap in the spectrum of $M_N$ pointed out by  Bai and Silverstein, it corresponds 
a gap in the spectrum of $A_NA_N^*$ which splits the spectrum of $A_NA_N^*$ exactly as that of $M_N$. We use the previous results  to characterize the outliers of spiked Information-Plus-Noise type models.
\end{abstract}

\section{Introduction}
Let $\{X_{ij}, i\in \mathbb{N}, j  \in \mathbb{N}\}$ be an infinite set of independent standardized complex random variables ($\mathbb{E}(X_{ij})=0, \mathbb{E}(\vert X_{ij}\vert^2 ) =1$) in some probability space. Define for any $\sigma>0$ and any nonnull integer numbers $n\leq N$,  the following matrix 
\begin{equation}\label{modele}M_N=( \sigma \frac{ X_N}{\sqrt{N}}+A_N)(\sigma \frac{ X_N}{\sqrt{N}}+A_N)^*\end{equation}
where $X_N=(X_{ij})_{1\leq i \leq n; 1\leq j \leq N}$ and $A_N$ is an $n\times N$ nonrandom matrix.
This model is referred to in the literature as the Information-Plus-Noise model.
For any Hermitian $n\times n$ matrix $Y$, denote by 
$$\lambda_1(Y) \geq \ldots \geq \lambda_n(Y)$$
\noindent the ordered eigenvalues of $Y$ and by $\mu_{Y}$  the empirical spectral measure of $Y$: $$\mu _{Y} := \frac{1}{n} \sum_{i=1}^n \delta_{\lambda _{i}(Y)}.$$  For a probability measure $\tau $ on $\R$, 
 denote by $g_\tau $ its Stieltjes transform defined for $z \in \C\setminus  \mbox{supp}(\tau)$ by 
$$g_\tau (z) = \int_\R \frac{d\tau (t)}{z-t}.$$
As $N$ tends to infinity, if $c_N=n/N \rightarrow c\in ]0,1]$, if the $X_{ij}$ are independent and identically distributed random variables and if the empirical spectral measure $\mu _{A_N A_N^*}$ of  $A_N A_N^*$ converges weakly  to some probability distribution $\nu \neq \delta_0$, Dozier and Silverstein established 
in \cite{DozierSilver} that almost surely the empirical spectral measure $\mu _{M_N}$ of $M_N$ converges weakly  towards a nonrandom distribution $\mu_{\sigma,\nu,c}$  which is characterized in terms of its Stieljes transform which satisfies the following equation:
 for any $z \in \mathbb{C}^+$,
\begin{equation}\label{TS} g_{\mu_{\sigma,\nu,c}}(z)=\int \frac{1}{(1-\sigma^2cg_{ \mu_{\sigma,\nu,c}}(z))z- \frac{ t}{1- \sigma^2 cg_{ \mu_{\sigma,\nu,c}}(z)} -\sigma^2 (1-c)}d\nu(t).\end{equation}
Note that, since for any $z \in \mathbb{C}^+$, $\Im ( z  g_{ \mu_{\sigma,\nu,c}}(z)) =-\int \frac{t \Im z}{\vert t-z\vert^2} d\mu_{\sigma,\nu,c}(z) \leq 0$ and $\Im (   g_{ \mu_{\sigma,\nu,c}}(z)) <0$, one can easily see that
the imaginary part of the denominator of the integrand in (\ref{TS}) is greater or equal to $\Im z$ so that the integral is well defined.
Note also that  in \cite{DozierSilver}, the authors proved that   the solution $m$ to the equation 
\begin{equation}m(z)=\int \frac{1}{(1-\sigma^2 c m(z))z- \frac{ t}{1- \sigma^2 cm(z)} -\sigma^2 (1-c)}d\nu(t), \mbox{for any $z \in \mathbb{C}^+$}\end{equation} 
\noindent is unique if $\Im m(z) <0$ and $\Im ( z  m(z))\leq 0$ (specifically if $m$ is the Stieljes transform of a probability measure supported on $[0;+\infty[$).\\
This result of convergence was extended to independent but non identically distributed  random variables by Xie in \cite{X}. (Note that, in \cite{HLN}, the authors investigated the case where $\sigma$ is replaced by  a bounded sequence of real numbers.)\\
In \cite{DozierSilver2}, Dozier and Silverstein investigated  the limiting spectral measure $\mu_{\sigma,\nu,c}$ and proved in particular that its distribution function is continuous (it has a continuous derivative on $\mathbb{R}\setminus\{0\}$ and no mass at zero).\\
  The support of the  probability measure $\mu_{\sigma,\mu _{A_N A_N^*},c_N}$ (that is the measure whose Stieljes transform satisfies  (\ref{TS}) where $\nu$ is replaced by $\mu _{A_N A_N^*}$ and $c$ by $c_N$)  plays a fundamental role in the study of the spectrum of $M_N$ (see \cite{BaiSilver, VLM, LV}).
Introducing,  for any probability measure $\tau$ on $[0,+\infty[$ and any $0<\gamma \leq 1$, $\rho>0$,  the function $\omega_{\rho,\tau,\gamma}$ defined in $\mathbb{R}\setminus \mbox{supp}(\mu_{\rho,\tau,\gamma})$
by 
 $$\omega_{\rho,\tau,\gamma}
(x )= 
x (1- \rho^2  \gamma g_{\mu_{\rho, \tau,\gamma}}(x))^2 -\rho^2 (1-\gamma)(1-\rho^2 \gamma g_{\mu_{\rho,\tau,\gamma}}(x)),$$
Bai and Silverstein established the following result.
\begin{theoreme}\label{pasde}\cite{BaiSilver}
Assume that 
\begin{enumerate}
\item $X_{ij}, i,j =1,2, \ldots$  are independent standardized random variables.
\item There exists a $K$ and a random variable $X$ with finite fourth moment for which there exists $x_0>0$ and integer number $n_0>0$ such that, for any $x >x_0$ and any integer numbers $n_1,n_2 >n_0$, we have
$$\frac{1}{n_1 n_2} \sum_{i\leq n_1, j\leq n_2}P\left( \vert X_{ij}\vert >x\right) \leq KP\left(\vert X \vert>x\right).$$
\item There exists a positive function $\Psi(x) \uparrow \infty$ as $ x \rightarrow \infty$ and $M>0$ such that 
$$\max_{ij} E\vert X_{ij}^2 \vert \Psi \left( \vert X_{ij}^2 \vert \right)\leq M.$$
\item $A_N$ is an $n\times N$ nonrandom matrix such that $\Vert A_N \Vert$ is uniformly bounded.
\item As $N$ tends to infinity, $c_N=n/N \rightarrow c\in ]0,1]$.
\item The empirical spectral measure $\mu _{A_N A_N^*}$ of  $A_N A_N^*$ converges weakly  to some probability distribution $\nu\neq \delta_0$.
\item
Let $[a,b] $ be such that the couple $ {\cal P}(\sigma)$ of the following properties is satisfied:
\begin{itemize}
\item (i) there exists $0<\delta<a$ such that for all large $N$, $ ]a-\delta; b+\delta[  \subset \mathbb{R}\setminus \rm{supp} (\mu_{\sigma,\mu _{A_N A_N^*},c_N})$,
 \item (ii)   $A_{Nj}$ denoting the matrix resulting from removing the $j$-th column from $A_N$, there exists $0<\tau<\delta$ and a positive $d<1$ such that for all $N$ large, the number of j's with no eigenvalues 
of $N/(N-1) A_{Nj}A_{Nj}^*$ appearing in   $ \omega_{{\sigma,\mu _{A_N A_N^*},c_N}}(]a -\tau, b +\tau[)$ is greater than $N-N^{d}$. 
\end{itemize}
\end{enumerate}
Then, $M_N$ being defined in (\ref{modele}),
 $$\mathbb P[\mbox{for all large N}, \rm{spect}(M_N) \subset  \R \setminus [a,b] ]=1.$$
\end{theoreme}
\begin{remarque} \label{incl}  Since $\mu_{\sigma,\mu _{A_N A_N^*},c_N}$ converges weakly towards $\mu_{\sigma,\nu,c}$ (this can be deduced from \cite{DozierSilver} and Theorem 1 in \cite{LV}), Assumption 7. (i) implies that 
$\forall 0< \tau< \delta$, $[a-\tau; b+\tau] \subset \mathbb{R} \setminus  \rm{supp}~ \mu_{\sigma,\nu,c}$.
\end{remarque}
Note that such a result was proved by a different approach in \cite{VLM} when the  $X_{ij}$ are independent gaussian variables without assuming condition (ii) of 7. in Theorem \ref{pasde}. Note that when the $X_{ij}$ are independent gaussian variables, it can be assumed that $A_N$ is such that 
\begin{equation}\label{diagonale} A_N=\begin{pmatrix}  a_1(N) ~~~~~~~~~~~~~~~~(0)\\
~~~(0)\\ ~~~~~~~~~~\ddots~~~~~~~~~~( 0)\\ ~(0)~~~~~~~~~~~~~~~~~~~~\\  ~~~~~~~~~~~~~a_{n}(N)~~~ ( 0  )    \end{pmatrix} \end{equation}  and, then, 
 condition (ii) of 7. in Theorem \ref{pasde} is not needed.\\

In  \cite{LV},  dealing with independent gaussian variables $X_{ij}$, P. Loubaton and P. Vallet established an exact separation phenomenon between the spectrum of $M_N$ and the spectrum of $A_NA_N^*$: 
 to a gap in the spectrum of $M_N$, it corresponds 
a gap in the spectrum of $A_NA_N^*$ which splits the spectrum of $A_N$ exactly as that of $M_N$. 
In this paper, we extend their result to the framework of non gaussian  Information-Plus-Noise type matrices investigated in \cite{BaiSilver} since we establish the following

\begin{theoreme}\label{sep} Assume conditions [1-7] of Theorem \ref{pasde} are satisfied.
If $c<1$, assume moreover that  $\omega_{\sigma,\nu,c}(b)>0$.
 Then for  $N$ large enough, $$\omega_{{\sigma,\nu,c}}([a,b])=[\omega_{{\sigma,\nu,c}}(a);\omega_{{\sigma,\nu,c}}(b)] \subset \mathbb{R} \setminus \mbox{supp}(\mu _{A_N A_N^*}).$$
With the convention that $\lambda _0(M_N)=\lambda _0(A_NA_N^*)=+\infty $ and 
$\lambda _{n+1}(M_N)=\lambda _{n+1}(A_NA_N^*)=-\infty $, for $N$ large enough, let  $i_N\in \{0,\ldots,n\}$
be such that
\begin{equation}\label{iN2}\lambda_{i_N+1}(A_N A_N^*) <\omega_{{\sigma,\nu,c}}(a) \mbox{~~ and ~~} \lambda_{i_N}(A_N A_N^*) > \omega_{{\sigma,\nu ,c}}(b).\end{equation}
Then
$$P[\mbox{for all large N}, \lambda_{i_N+1}(M_N) <a\mbox{~and} ~ \lambda_{i_N}(M_N)> 
b] = 1.$$
\end{theoreme}

Such an exact separation phenomenon was previously exhibited  for sample covariance matrices in \cite{BaiSil99} and later proved for deformed Wigner matrices in \cite{CDFF}.
The technics of \cite{LV} completely differ from the approach used in the present paper which uses an argument similar to \cite{CDFF} which consists in  introducing a sequence of matrices interpolating between $A_NA_N^*$ and $M_N$.
Before using this approach, we carry on with the study of the  support of the limiting spectral measure previously investigated in \cite{DozierSilver2} and later
in  \cite{VLM, LV}. The main  results about  the limiting spectral measure can be summarized in the following theorem.

\begin{theoreme}\label{caractfinale}
Let $\gamma $ be in $]0;1]$, $\rho$ be in $]0; +\infty[$ and $\tau$ be a compactly supported probability measure on $[0;+\infty[$. 
Define differentiable functions   $\omega_{\rho,\tau,\gamma}$ and $\Phi_{\rho,\tau,\gamma}$  on respectively $ \mathbb{R}\setminus \mbox{supp}(\mu_{\rho,\tau,\gamma})$ and $ \mathbb{R}\setminus \mbox{supp}(\tau)$ by setting  $$\omega_{\rho,\tau,\gamma} :\begin{array}{ll} \mathbb{R}\setminus \mbox{supp}(\mu_{\rho,\tau,\gamma}) \rightarrow \mathbb{R}\\
x \mapsto 
x (1- \rho^2 \gamma g_{ \mu_{\rho,\tau,\gamma}}(x))^2 -\rho^2 (1-\gamma)(1-\rho^2 \gamma g_{\mu_{\rho,\tau,\gamma}}(x))\end{array}$$
and  $$\Phi_{\rho,\tau,\gamma} :\begin{array}{ll} \mathbb{R}\setminus \mbox{supp}(\tau) \rightarrow \mathbb{R}\\
x \mapsto 
x (1+\gamma \rho^2g_{ \tau}(x))^2 + \rho^2 (1-\gamma) (1+ \gamma \rho^2 g_\tau(x))\end{array}.$$
Set $$ {\cal E}_{\rho,\tau,\gamma}:=\left\{ x \in  \mathbb{R}\setminus \mbox{supp}(\tau),  \Phi_{\rho,\tau,\gamma}^{'}(x) >0, g_\tau(x) >-\frac{1}{\rho^2\gamma}\right\}.$$
We have the following results.
 \begin{itemize} 
\item[A)] \begin{enumerate}
 \item If $\gamma<1$ then $0 \notin \rm{supp}(\mu_{\rho,\tau,\gamma}).$
\item  We have  $0 \notin \rm{supp}(\mu_{\rho,\tau,1}) $ if and only if $0\notin \rm{supp}( \tau)$ and  $g_\tau(0) >-\frac{1}{\rho^2}$.
\end{enumerate}
\item[B)] 
$\omega_{\rho,\tau,\gamma}$ is an increasing  analytic diffeomorphism with positive derivative  from $\mathbb{R}\setminus \mbox{supp}(\mu_{\rho,\tau,\gamma})$ to ${\cal E}_{\rho,\tau,\gamma}$, with inverse $\Phi_{\rho,\tau,\gamma}$.
\item[C)]  Let $a$ be in $\mathbb{R} \setminus \mbox{supp}(\mu_{\rho,\tau,\gamma}) $ such that $\omega_{\rho,\tau,\gamma}(a) \leq 0.$ Then
$a$ is to the  left  of $\mbox{supp}(\mu_{\rho,\tau,\gamma})$.

\item[D)] Assume that the support of $\tau$ is a finite union of disjoint (possibly degenerate) closed bounded intervals.
There exists a nonnul integer number  $p$ and  $u_1< v_1<u_2<\ldots <u_p < v_p$  (depending on $\rho,\tau,\gamma$) such that 
$${{\cal E}_{\rho,\tau,\gamma}}=]-\infty;u_1[\cup_{l=1}^{p-1}]v_l;u_{l+1}[\cup 
]v_p;+\infty[.$$
We have $$\mbox{supp}(\tau) \subset \cup_{l=1}^{p} [u_l;v_l]$$
and for  each $l \in \{1,\ldots,p\}$,  $[u_l;v_l]\cap\mbox{supp}(\tau) \neq \emptyset$.

\noindent Moreover,
$$
 \mbox{supp}(\mu_{\rho,\tau,\gamma})
=\cup_{l=1}^p [\Phi_{\rho,\tau,\gamma}(u_l^-);\Phi_{\rho,\tau,\gamma}(v_l^+)],$$
with \\

$\Phi_{\rho,\tau,\gamma}(u_1^-) < \Phi_{\rho,\tau,\gamma}(v_1^+)< \Phi_{\rho,\tau,\gamma}(u_2^-) < \Phi_{\rho,\tau,\gamma}(v_2^+)$
$$\hspace*{6cm}
<\cdots< \Phi_{\rho,\tau,\gamma}(u_p^-) < \Phi_{\rho,\tau,\gamma}(v_p^+),$$
where $ \Phi_{\rho,\tau,\gamma}(u_l^-) =\lim_{u\uparrow u_l} \Phi_{\rho,\tau,\gamma}(u)$ and $ \Phi_{\rho,\tau,\gamma}(v_l^+) =\lim_{u\downarrow v_l} \Phi_{\rho,\tau,\gamma}(u)$.
Finally, for  each $l \in \{1,\ldots,p\}$, \begin{equation}\label{palier}\mu _{\rho, \tau, \gamma }([\Phi _{\rho , \tau,\gamma }(u_l^-), \Phi_{\rho, \tau,\gamma }(v_l^+)]) = \tau ([u_l, v_l]).\end{equation}

\end{itemize}
\end{theoreme}

\begin{remarque}\label{FW}
The restriction in Theorem \ref{sep}, ``if $c<1$, $\omega_{\sigma,\nu,c}(b)>0$", could be relaxed if  Theorem \ref{pasde}  of \cite{BaiSilver}
could ever be extended to the case $a=0$.
 We will prove  in Theorem \ref{caractfinale}  that, if $c<1$, the minimum of the support of $\mu_{\sigma,\nu,c}$ is
 positive. Let us denote by $x_0$
this minimum. We will see that  if $b\in [0;+\infty[\setminus \rm{supp} (\mu_{\sigma,\nu,c}) $  does not belong to  $[0;x_0[$, then $\omega_{\sigma,\nu,c}(b)>0$ (see Lemma \ref{CSomega}). 
Nevertheless this is not a necessary condition for $\omega_{\sigma,\nu,c}(b)>0$ to hold as the following example shows. Let $c<1$ and $\sigma>0$.
Let $w>0$ and define $d\nu(t)=\frac{3}{(y-w)^3} \1_{[w;y]}(t) (w-t)^2 dt$ where $y>w$ is choosen such that $L=y-w$ is large enough to satisfy 
$$\lim_{x \uparrow w} -g_\nu(x)=\frac{3}{2L} < \frac{1}{\sigma^2c}$$
and $$\lim_{x \uparrow w} \Phi^{'}_{\sigma,\nu,c}(x)=-\left\{2w (1-\sigma^2c \frac{3}{2L}) +\sigma^2 (1-c) \right\} \frac{3 \sigma^2 c}{L^2} + (1-\sigma^2c \frac{3}{2L})^2 >0.$$
Then, with the notations of D) in Theorem \ref{caractfinale}, we have $u_1=w$ and the interval  $]\Phi_{\sigma,\nu,c}(0); \Phi_{\sigma,\nu,c}(w^{-})[$ is to the left of 
 the support of $\mu _{\sigma,\nu,c}$ whereas its image by $\omega_{\sigma,\nu,c}$ is $]0;w[$.

\end{remarque}


The previous Theorem \ref{sep} and Theorem \ref{caractfinale} allow us to investigate the spectrum of  spiked 
models when the perturbation matrix can be written as in (\ref{diagonale}) and to obtain in Theorem \ref{ThmASCV2} a  description of the convergence of the
eigenvalues of $M_N$ depending on the location of the spikes of the perturbation  with respect to $ {\cal E}_{\sigma,\nu,c}$
and to the connected components of the support of $\nu$. This extends  previous results in \cite{LV} and \cite{BR} involving the  Gaussian case and finite rank perturbations.

This paper is organized as follows.
 Section \ref{lsm} is devoted to the proof of Theorem \ref{caractfinale}. Theorem \ref{sep} is proved in Section \ref{thsep}. Section \ref{spi}  investigates the spectrum of  spiked 
models. 
\section{The  support of the limiting spectral measure}\label{lsm}

In this section, we split  the results of  Theorem \ref{caractfinale} into different propositions that we  prove successively.\\
Firstly, in subsection \ref{ds} we  are going to take up the arguments and results of \cite{DozierSilver2}
 that we will develop here for the reader's convenience in order to state, in Proposition \ref{inverserec},  the characterization of the  complement in $\mathbb{R}\setminus \{0\}$ of the support of $\mu_{\rho,\tau,\gamma}$. \\Secondly, in subsection \ref{CNS0}, , we establish a necessary and sufficient condition for the inclusion of zero in the support of $\mu_{\rho,\tau,\gamma}$, in  Propositions \ref{zeroinclus1} and   \ref{zeroinclus}, which allows us to put foward a
complete characterization of the  complement of the support of $\mu_{\rho,\tau,\gamma}$  in $\mathbb{R}$,  in Proposition \ref{supportct}.\\
We prove the strict monoticity of  $\Phi_{\rho,\tau,\gamma}$ 
 on the set ${\cal E}_{\rho,\tau,\gamma}$ in subsection \ref{1}. Note that we first prove it for $\gamma=1$  in Proposition \ref{phicroitglobalement1} and then use 
   a relationship between  the  complements of the supports of $\mu_{\rho,\tau,\gamma}$
and   $\mu_{\rho \sqrt{\gamma},\tau,1}$, established in  Proposition \ref{passagede1ac}, to extend the result to any $\gamma<1$
 in Proposition \ref{phicroitglobalement}.\\
In subsection \ref{conc}, we explain how to deduce Theorem \ref{caractfinale} from the previous results.

\subsection{Characterization of the  complement in $\mathbb{R}\setminus \{0\}$ of the support of $\mu_{\rho,\tau,\gamma}$}\label{ds}


Actually it is possible to deduce from results of \cite{DozierSilver2} the following characterization of the complement of $\mbox{supp}(\mu_{\rho,\tau,\gamma})\cup \{0\}$ and the following relationships between $\Phi_{\rho,\tau,\gamma}$ and $\omega_{\rho,\tau,\gamma}$ .

\begin{proposition}\label{inverserec}
~~
\begin{itemize}
\item(i) ~~$\forall x \in   \mathbb{R}\setminus \{ \mbox{supp}(\mu_{\rho,\tau,\gamma})\cup\{0\}\}$, we have
$ \omega_{\rho,\tau,\gamma}(x) \in  \mathbb{R}\setminus \mbox{supp}(\tau), \omega_{\rho,\tau,\gamma}^{'}(x)>~0, \\  ~~\Phi_{\rho,\tau,\gamma}^{'}( \omega_{\rho,\tau,\gamma}(x))>0, $~~$g_\tau( \omega_{\rho,\tau,\gamma}(x)) >-\frac{1}{\rho^2\gamma}$ and 
\begin{equation}\label{conjdsunsens}\Phi_{\rho,\tau,\gamma}(\omega_{\rho,\tau,\gamma}(x))=x.\end{equation}
\item(ii) $ ~~\forall x \in   \mathbb{R}\setminus \mbox{supp}(\tau)  \mbox{  such that }\Phi_{\rho,\tau,\gamma}^{'}(x) >0  \mbox{ and } g_\tau(x) >-\frac{1}{\rho^2\gamma},$ we have  \\ $\Phi_{\rho,\tau,\gamma}(x) \in  \mathbb{R}\setminus  \mbox{supp}(\mu_{\rho,\tau,\gamma}) \mbox{~~and~~}$
 \begin{equation}\label{omegacircphi}\omega_{\rho,\tau,\gamma}\left(\Phi_{\rho,\tau,\gamma}(x)\right)=x.\end{equation}
\item(iii) Therefore \\

$ \mathbb{R}\setminus \{ \mbox{supp}(\mu_{\rho,\tau,\gamma})\cup\{0\}\}  =$
\begin{equation}\label{support} \Phi_{\rho,\tau,\gamma}\left\{ x \in   \mathbb{R}\setminus \mbox{supp}(\tau), \Phi_{\rho,\tau,\gamma}(x) \neq 0, \Phi_{\rho,\tau,\gamma}^{'}(x) >0, g_\tau(x) >-\frac{1}{\rho^2\gamma}\right\}. \end{equation}
\end{itemize}
\end{proposition}

\begin{proof}
\noindent  (\ref{TS}) may be rewritten,  for any $z \in \mathbb{C}^+$,
\begin{equation}\label{TS2}\frac{ g_{ \mu_{\rho,\tau,\gamma}}(z)}{1- \rho^2 \gamma g_{ \mu_{\rho,\tau,\gamma}}(z)}=g_\tau\left[ z(1- \rho^2 \gamma g_{ \mu_{\rho,\tau,\gamma}}(z))^2 -\rho^2 (1-\gamma) (1-\rho^2 \gamma g_{\mu_{\rho,\tau,\gamma}}(z))\right]\end{equation}
According to Theorem 2.1 and Lemma 2.1   in \cite{DozierSilver2}, for any  $ x $ in $ \mathbb{R}\setminus\{0\}$, one has $\lim_{z\in \mathbb{C}^+\rightarrow  x }g_{\mu_{\rho,\tau,\gamma}}(z):=g_{\mu_{\rho,\tau,\gamma}}(x)$ exists and
\begin{equation}\label{Re} \forall z \in \mathbb{C}^+ \cup \mathbb{R}\setminus\{0\}, ~   \Re (\frac{1}{\rho^2 \gamma}-g_{\mu_{\rho,\tau,\gamma}}(z))>0.\end{equation} Moreover,  according to Theorem 3.2  in \cite{DozierSilver2}, 
\begin{equation}\label{imageomega} \mbox{if  } x \in  \mathbb{R}\setminus \{ \mbox{supp}(\mu_{\rho,\tau,\gamma})\cup\{0\}\} \mbox{ then~} \omega_{\rho,\tau,\gamma}(x) \in  \mathbb{R}\setminus \mbox{supp}(\tau)
. \end{equation}
Therefore it makes sense to let $z \in  \mathbb{C}^+ $ tend to $x \in  \mathbb{R}\setminus \{\mbox{supp}(\mu_{\rho,\tau,\gamma}) \cup \{0\} \}$ in (\ref{TS2}) and get
  for any $x \in  \mathbb{R}\setminus \{\mbox{supp}(\mu_{\rho,\tau,\gamma}) \cup \{0\} \}$,
\begin{equation}\label{TS3}\frac{ g_{ \mu_{\rho,\tau,\gamma}}(x)}{1- \rho^2 \gamma g_{ \mu_{\rho,\tau,\gamma}}(x)}=g_\tau\left(\omega_{\rho,\tau,\gamma}(x)\right).\end{equation}
Multiplying both sides of (\ref{TS3}) by $\rho^2 \gamma$ and adding 1, it readily comes that 

\begin{equation}\label{TS3bis} \frac{1}{1-\rho^2 \gamma g_{ \mu_{\rho,\tau,\gamma}}(x)}=1 +\rho^2 \gamma g_\tau\left(\omega_{\rho,\tau,\gamma}(x)\right)\end{equation}
and then \begin{equation} \label{TS4} 1- \rho^2 \gamma g_{ \mu_{\rho,\tau,\gamma}}(x)= \frac{1}{1 + \rho^2 \gamma g_\tau\left(\omega_{\rho,\tau,\gamma}(x)\right)}. \end{equation}
Replacing $ 1- \rho^2 \gamma g_{ \mu_{\rho,\tau,\gamma}}(x)$ in the expression of $\omega_{\rho,\tau,\gamma}(x)$ by the right hand side of (\ref{TS4}),
it follows that  $$\omega_{\rho,\tau,\gamma}(x)= \frac{x}{\left(1 +\rho^2 \gamma g_\tau\left(\omega_{\rho,\tau,\gamma}(x)\right)\right)^2} - \frac{\rho^2(1-\gamma)}{\left(1 +\rho^2 \gamma g_\tau\left(\omega_{\rho,\tau,\gamma}(x)\right)\right)}$$
and finally that  \begin{equation}\label{Inverse}\Phi_{\rho,\tau,\gamma}(\omega_{\rho,\tau,\gamma}(x))=x\end{equation}
Hence any $x \in  \mathbb{R}\setminus \{\mbox{supp}(\mu_{\rho,\tau,\gamma}) \cup \{0\} \}$ can be written as $\Phi_{\rho,\tau,\gamma}(u)$ where $u =\omega_{\rho,\tau,\gamma}(x) \in 	\mathbb{R}\setminus \mbox{supp}(\tau) $.
Let us prove that $ \omega_{\rho,\tau,\gamma}^{'}(x)>~0$, $\Phi_{\rho,\tau,\gamma}^{'}(\omega_{\rho,\tau,\gamma}(x))>0$ and $g_\tau (\omega_{\rho,\tau,\gamma}(x))>-\frac{1}{\rho^2 \gamma}$.  By differentiating both sides of (\ref{TS3bis}) we obtain that  for any 
$x \in  \mathbb{R}\setminus \{\mbox{supp}(\mu_{\rho,\tau,\gamma}) \cup \{0\} \}$  \begin{equation}\frac{ g^{'}_{ \mu_{\rho,\tau,\gamma}}(x)}{(1- \rho^2  \gamma g_{ \mu_{\rho,\tau,\gamma}}(x))^2}=\omega_{\rho,\tau,\gamma}^{'} (x) g_\tau^{'} \left(\omega_{\rho,\tau,\gamma}(x)\right).\end{equation}
Therefore since $g^{'}_{ \mu_{\rho,\tau,\gamma}}(x)<0$ and $g_\tau^{'} \left(\omega_{\rho,\tau,\gamma}(x)\right)<0$ we can deduce that 
\begin{equation} \label{croissanceomega} \mbox{for any 
$x \in  \mathbb{R}\setminus \{\mbox{supp}(\mu_{\rho,\tau,\gamma}) \cup \{0\} \}$} ,
\omega_{\rho,\tau,\gamma}^{'} (x) >0. \end{equation}
Now by differentiating both sides of (\ref{Inverse}) we obtain that 
\begin{equation} \label{croissancePhi} \mbox{for any 
$x \in  \mathbb{R}\setminus \{\mbox{supp}(\mu_{\rho,\tau,\gamma}) \cup \{0\} \}$} ,
\Phi_{\rho,\tau,\gamma}^{'}(\omega_{\rho,\tau,\gamma}(x))>0. \end{equation}
Let us prove now that
\begin{equation} \label{mingnu} \mbox{ for  $x \in  \mathbb{R}\setminus \{\mbox{supp}(\mu_{\rho,\tau,\gamma}) \cup \{0\} \}$},  g_\tau (\omega_{\rho,\tau,\gamma}(x))>-\frac{1}{\rho^2 \gamma}.\end{equation}
According to Lemma 2.1 (c) in \cite{DozierSilver2} ,  $g_{\mu_{\rho,\tau,\gamma}}(x)<\frac{1}{\rho^2 \gamma}$.  Moreover $ y\mapsto \frac{y}{1-\rho^2 \gamma y} $ is 
 increasing from $]-\infty ; \frac{1}{\rho^2 \gamma}[$ onto $]-\frac{1}{\rho^2 \gamma}; + \infty[$.
The result readily follows using (\ref{TS3}). 
%
The proof of (i) is complete.\\

Now let $x$ be in  $ \mathbb{R}\setminus \mbox{supp}(\tau)$, such that $  \Phi_{\rho,\tau,\gamma}^{'}(x) >0$ and $ g_\tau(x) >-\frac{1}{\rho^2 \gamma}$.
Following \cite{DozierSilver2}, let us prove that $\Phi_{\rho,\tau,\gamma}(x)$ belongs to $ \mathbb{R}\setminus \mbox{supp}(\mu_{\rho,\tau,\gamma})$ and  $ \omega_{\rho,\tau,\gamma}(\Phi_{\rho,\tau,\gamma}(x))=x$. Let $]l_1;l_2[ \subset [L_1;L_2] \subset  \mathbb{R}\setminus \mbox{supp}(\tau)$
such that $x \in ]l_1;l_2[ $ and for any $v \in   ]l_1;l_2[$,  $\Phi_{\rho,\tau,\gamma}^{'}(v) >0, g_\tau(v) >-\frac{1}{\rho^2 \gamma}.$
$ g_\tau$ is 
decreasing on $  ]l_1;l_2[$ and maps $  ]l_1;l_2[$ onto some interval $]d_1;d_2[ \subset ]-\frac{1}{\rho^2 \gamma};+\infty[.$
$h: b \mapsto \frac{1}{\rho^2 \gamma} \left(\frac{1}{b} -1\right) $ is a 
decreasing function from $]0; + \infty [$ onto $]-\frac{1}{\rho^2 \gamma}; + \infty[$.
Hence there is an interval $ ]k_1;k_2[ \subset ]0;+ \infty[$ such that $h$ is a one-to-one correspondence from $ ]k_1;k_2[$ to  $]d_1;d_2[ $. Therefore 
$ g_\tau^{-1}\circ h $ is is a one-to-one correspondence from $ ]k_1;k_2[$ to  $]l_1;l_2[ $.
For any $v$ in $ ]l_1;l_2[$, there exists a unique  $k$ in $ ]k_1;k_2[$ such that $v= g_\tau^{-1}(h(k))$; then $k=k(v)=h^{-1}(g_\tau(v))=\frac{1}{1+\rho^2 \gamma g_\tau(v)}$.
Moreover  \begin{eqnarray} p(k)&:= &\frac{1}{k^2}  g_\tau^{-1}( h(k))+ \frac{1}{k} \rho^2 (1-\gamma) \nonumber \\& = &(1+ \rho^2 \gamma g_\tau (v))^2 v +
(1+\rho^2 \gamma g_\tau(v)) \rho^2 (1-\gamma) \nonumber\\& =&\Phi_{\rho,\tau,\gamma}(v) \label{kv}.\end{eqnarray}
Differentiating both sides of (\ref{kv}) 
 in $v$ we obtain that for any $v \in ]l_1;l_2[$,  $ p^{'}(k(v))k^{'}(v)=\Phi_{\rho,\tau,\gamma}^{'}(v)$ with $k^{'}(v)=- \rho^2 \gamma \frac{g_\tau^{'} (v)}{ (1+\rho^2 \gamma g_\tau (v))^2}>0$ .
Therefore $\Phi_{\rho,\tau,\gamma}^{'}(v) >0$ implies $ p^{'}(k(v))>0$. According to Theorem 3.3 in \cite{DozierSilver2}, we can conclude that 
$ p(k(x))= \Phi_{\rho,\tau,\gamma}(x) \in \mathbb{R} \setminus  \mbox{supp}(\mu_{\rho,\tau,\gamma})$ and \begin{equation}\label{kp}k=1-\rho^2 \gamma g_{\mu_{\rho,\tau,\gamma}}(p(k)).\end{equation} 
Moreover (\ref{kp}) implies
$$\frac{1}{1+\rho^2 \gamma g_\tau(x)}=1-\rho^2 \gamma g_{\mu_{\rho,\tau,\gamma}}(\Phi_{\rho,\tau,\gamma}(x)).$$  It readily follows that $$\Phi_{\rho,\tau,\gamma}(x)(1-\rho^2 \gamma g_{\mu_{\rho,\tau,\gamma}}(\Phi_{\rho,\tau,\gamma}(x))^2- \rho^2(1-\gamma)(1-\rho^2 \gamma g_{\mu_{\rho,\tau,\gamma}}(\Phi_{\rho,\tau,\gamma}(x)))=x$$ that is $\omega_{\rho,\tau,\gamma}(\Phi_{\rho,\tau,\gamma}(x))=x$.
The proof of (ii) is complete.\\
(iii) is a straightforward consequence of (i) and (ii). \end{proof}
\subsection{A necessary and sufficient condition for the inclusion of zero in the support of $\mu_{\rho,\tau,\gamma}, 0<  \gamma \leq 1$, and  a complete characterization of the  complement of the support }\label{CNS0}
When $\gamma=1$, we have the following equivalence.
\begin{proposition}\label{zeroinclus1}
$$0 \notin \rm{supp}(\mu_{\rho,\tau,1}) \mbox{~if and only if ~}   0\notin \rm{supp}( \tau) \mbox{~and~} g_\tau(0) >-\frac{1}{\rho^2}.$$
\end{proposition}
\begin{proof}~~

\begin{itemize}
\item[$(\Longleftarrow)$] If  $0\notin \rm{supp}( \tau)$ and if $g_\tau(0) >-\frac{1}{\rho^2}$ then $\Phi^{'}_{\rho,\tau,1}(0)=(1+ \rho^2 g_\tau(0))^2 >0$. According to Proposition \ref{inverserec} (ii), we have 
$\Phi_{\rho,\tau,1}(0) (=0) \in \mathbb{R} \setminus  \rm{supp}(\mu_{\rho,\tau,1}).$ 
\item[$(\Longrightarrow)$] Assume that $0 \notin \rm{supp}(\mu_{\rho,\tau,1})$. Let $\epsilon>0$ be such that  $]-\epsilon; +\epsilon[ \subset \mathbb{R} \setminus   \rm{supp}(\mu_{\rho,\tau,1})$. According to Proposition  \ref{inverserec} (i), we have for $x \neq 0$ in $]-\epsilon;\epsilon[$,  $\omega_{\rho,\tau,1}(x)
\in \mathbb{R} \setminus \rm{supp}( \tau)$, $\omega^{'}_{\rho,\tau,1}(x)>0$, 
\begin{equation}\label{minoration} g_\tau(\omega_{\rho,\tau,1}(x))> -\rho^{-2}, \end{equation}
and $\Phi_{\rho,\tau,1}(\omega_{\rho,\tau,1}(x) ) =x$. Note that $
\omega_{\rho,\tau,1}(0)=0$.
Assume first that  $0\in   \rm{supp}( \tau)$. Since $]\omega_{\rho,\tau,1} (-\epsilon); 0[\cup ]0; \omega_{\rho,\tau,1} (\epsilon)[ \subset \mathbb{R} \setminus  \rm{supp}( \tau)$, we can conclude that $\tau$ has a mass at zero contradicting (\ref{minoration}).
 Therefore $0\notin  \rm{supp}( \tau)$. We have  $\Phi_{\rho,\tau,1}(\omega_{\rho,\tau,1}(0) ) =0$.
Thus,  we have for any $x$ in $]-\epsilon; +\epsilon[ $,
\begin{equation}\label{compositionid}
\Phi_{\rho,\tau,1}(\omega_{\rho,\tau,1}(x) ) =x.
\end{equation}
Differentiating  (\ref{compositionid}) yields $\Phi^{'}_{\rho,\tau,1}(0) =(1+ \rho^2 g_\tau(0))^2\neq 0$ and letting $x$ goes to zero in (\ref{minoration}) leads to $ g_\tau(0)\geq -\rho^{-2}$. Therefore, we have $ g_\tau(0)> -\rho^{-2}$. 
\end{itemize}
\end{proof}
When $0< \gamma< 1$, it turns out that the minimum of the support of $\mu_{\rho,\tau,\gamma}$ is positive as the following Proposition \ref{zeroinclus} shows.
The next straightforward  lemma will be used in the proof of Proposition \ref{zeroinclus}.
\begin{lemme}\label{annulation}
If $x \in \mathbb{R} \setminus  \rm{supp}( \tau)$ is such that $\Phi_{\rho,\tau,\gamma}^{'}(x)=0$ and $1 + \gamma \rho^2 g_\tau(x) >0$ then $\Phi_{\rho,\tau,\gamma}(x)>0$.
\end{lemme}
\begin{proof} Note that 
\begin{eqnarray*}
\left(1 + \gamma \rho^2 g_\tau(x) \right) \Phi_{\rho,\tau,\gamma}^{'}(x)&=& \left(1 + \gamma \rho^2 g_\tau(x) \right)^3 + 2 \Phi_{\rho,\tau,\gamma}(x) \gamma \rho^2 g_\tau^{'}(x) \\&& -\rho^4 (1-\gamma)\gamma \left(1 + \gamma \rho^2 g_\tau(x) \right) g_\tau^{'}(x).
\end{eqnarray*}
Since $g_\tau^{'}(x)<0$, it is then clear that, if $ \Phi_{\rho,\tau,\gamma}(x) \leq 0$ and $1 + \gamma \rho^2 g_\tau(x) >0$ then $\Phi^{'}_{\rho,\tau,\gamma}(x)>0$.
The lemma follows. \end{proof}
\begin{proposition}\label{zeroinclus}  If $0<\gamma<1$ then $0 \notin \rm{supp}(\mu_{\rho,\tau,\gamma}).$
\end{proposition}
\begin{proof} 
 Since $\lim_{x\rightarrow - \infty} \Phi_{\rho,\tau,\gamma}^{'}(x)= 1$ and  $\lim_{x\rightarrow - \infty} g_\tau(x)= 0$,  we can define  $$x_0= \sup\left\{x^{'} \leq \min \rm{supp}( \tau);
\forall x < x^{'},   g_\tau(x) >-\frac{1}{\rho^2 \gamma}, \Phi^{'}_{\rho,\tau,\gamma}(x)>0\right\}.$$
Since moreover, according to Proposition \ref{inverserec} (ii), we have $$]-\infty; \Phi_{\rho,\tau,\gamma}(x_0^{-})[ \subset \mathbb{R} \setminus  \rm{supp}(\mu_{\rho,\tau,\gamma}),$$ we are going to prove that $ \Phi_{\rho,\tau,\gamma}(x_0^{-}) >0.$
\begin{itemize}
\item Assume that $\lim_{x \uparrow \min \rm{supp}( \tau)} g_\tau (x)\geq - \frac{1}{\rho^2 \gamma}$. Therefore we have for any $x < \min  \rm{supp}( \tau)$, $g_\tau(x) > - \frac{1}{\rho^2 \gamma}.$ 
\begin{enumerate} \item If $x_0 <\min  \rm{supp}( \tau)$, we can conclude that  $\Phi_{\rho,\tau,\gamma}^{'}(x_0)=0$ and $ \left(1 + \gamma \rho^2 g_\tau(x_0) \right)>0$.
Then Lemma \ref{annulation} yields $\Phi_{\rho,\tau,\gamma}(x_0)>0$. \\
\item Assume that   $x_0 =\min  \rm{supp}( \tau)$ and that $\lim_{x \uparrow \min \rm{supp}( \tau)} g_\tau (x)= - \frac{1}{\rho^2 \gamma}$. 
Let $x_2$ be in $]-\infty; \min \rm{supp}( \tau)[$ and $y_2=g^{'}_\tau (x_2) <0.$ 
Let $x_1$ be such that $x_2<x_1 < \min \rm{supp}( \tau)$,
$\vert 2x_1 \left(1 + \gamma \rho^2 g_\tau(x_1) \right)\vert  < \frac{\rho^2 (1-\gamma)}{2}$ and 
$   \left(1 + \gamma \rho^2 g_\tau(x_1) \right)^2 < -y_2 \frac{\rho^4(1-\gamma)\gamma}{4}.$
Note that for any $ x < \min \rm{supp}( \tau)$, $g^{''}_\tau (x) =2 \int \frac{1}{(x-t)^3} d\tau(t) \leq 0$ and thus $g^{'}_\tau (x_1)\leq y_2.$ 
We have 
\begin{eqnarray*}
 \Phi_{\rho,\tau,\gamma}^{'}(x_1)&=& \left[2x_1 \left(1 + \gamma \rho^2 g_\tau(x_1) \right) + \rho^2 (1-\gamma) \right] \rho^2 \gamma   g_\tau^{'}(x_1)\\&& +\left(1 + \gamma \rho^2 g_\tau(x_1) \right) ^2 \\
&<&  \frac{\rho^4(1-\gamma)\gamma}{2}g_\tau^{'}(x_1) + \left(1 + \gamma\rho^2 g_\tau(x_1) \right)^2\\&\leq& \frac{\rho^4(1-\gamma)\gamma}{2} y_2  -y_2 \frac{\rho^4(1-\gamma)\gamma}{4}=  \frac{\rho^4(1-\gamma)\gamma}{4} y_2 <0,
\end{eqnarray*}
which leads to a contradiction.
\item Assume that   $x_0 =\min  \rm{supp}( \tau)$ and that $\lim_{x \uparrow \min \rm{supp}( \tau)} g_\tau (x)> - \frac{1}{\rho^2 \gamma}$. 
Then, it is clear by the very definition of $\Phi_{\rho,\tau,\gamma}$ that $ \Phi_{\rho,\tau,\gamma}(x_0^{-}) >0.$
\end{enumerate}

\item  If $\lim_{x \uparrow \min \rm{supp}( \tau)} g_\tau (x) < - \frac{1}{\rho^2 \gamma}$,
 there exists $x_1< \min  \mbox{supp}(\tau)$ such that  $g_\tau (x_1) = - \frac{1}{\rho^2 \gamma}$ and thus $\Phi_{\rho,\tau,\gamma}^{'}(x_1)=
\gamma(1-\gamma)\rho^4 g^{'}_\tau (x_1)<0.$ We can conclude that $x_0< x_1$, $\Phi_{\rho,\tau,\gamma}^{'}(x_0)=0$ and $ \left(1 + \gamma \rho^2 g_\tau(x_0) \right)>0$.
Then Lemma \ref{annulation} yields $\Phi_{\rho,\tau,\gamma}(x_0)>0$. 
\end{itemize}
\end{proof}
\begin{remarque} We presented  a proof of Proposition \ref{zeroinclus} which did not require non elementary background. Nevertheless, there is a shorter proof using free probability theory as follows.
Choosing gaussian entries for $X_N$, it is easy  to see that $\mu_{\rho,\tau,\gamma}$ is the distribution of a Free Wishart process $FW(\frac{1}{\gamma}; \tau)$,
introduced in \cite{CDM}, at time $t=\rho^2\gamma$ (see Section 2.2 p 421 in \cite{CDM}).
 It is proved in Proposition 2.2 in \cite{CDM} that if $\gamma<1$, the minimum of the support of $\mu_{\rho,\tau,\gamma}$  is strictly positive; therefore 
 $0 \notin \rm{supp}(\mu_{\rho,\tau,\gamma}).$
\end{remarque}
We are now in position to establish the following proposition.
\begin{proposition}\label{supportct} For any $0<\gamma\leq 1$,
$$\mathbb{R}\setminus \mbox{supp}(\mu_{\rho,\tau,\gamma}) = \Phi_{\rho,\tau,\gamma}\left\{ x \in   \mathbb{R}\setminus \mbox{supp}(\tau) , \Phi_{\rho,\tau,\gamma}^{'}(x) >0, g_\tau(x) >-\frac{1}{\rho^2 \gamma}\right\}.$$
\end{proposition}
\begin{proof}
When $\gamma=1$, Proposition \ref{supportct} readily follows from 
Proposition \ref{inverserec} and Proposition \ref{zeroinclus1}. Therefore we focus on the case $\gamma<1$.
According to Proposition \ref{zeroinclus}, $0\notin \mbox{supp}(\mu_{\rho,\tau,\gamma})$.  Since $x \mapsto \frac{1}{\rho^2 \gamma}-g_{\mu_{\rho,\tau,\gamma}}(x)$ is strictly increasing on $]-\infty; \min \mbox{supp}(\mu_{\rho,\tau,\gamma})[$
and according to (\ref{Re}), $ \forall x \in \mathbb{R}\setminus\{0\}, ~   \Re (\frac{1}{\rho^2 \gamma}-g_{\mu_{\rho,\tau,\gamma}}(x))>0$, we can readily deduce that  
$\frac{1}{\rho^2 \gamma}-g_{\mu_{\rho,\tau,\gamma}}(0)>0$ and then \begin{equation}\label{enzero}\omega_{\rho,\tau,\gamma}(0)<0.\end{equation}
Let $\epsilon>0$ be such that $[-\epsilon; \epsilon]\subset \mathbb{R}\setminus \mbox{supp}(\mu_{\rho,\tau,\gamma})$ and $\omega_{\rho,\tau,\gamma}(\epsilon)<0$.
According to Proposition \ref{inverserec} (i), 
 we have for any $ x \in [-\epsilon;\epsilon]\setminus\{0\},$ 
$ \omega_{\rho,\tau,\gamma}(x) \in  \mathbb{R}\setminus \mbox{supp}(\tau), \omega_{\rho,\tau,\gamma}^{'}(x)>~0, $ $\Phi_{\rho,\tau,\gamma}^{'}( \omega_{\rho,\tau,\gamma}(x))>0, $~~$g_\tau( \omega_{\rho,\tau,\gamma}(x)) >-\frac{1}{\rho^2\gamma}$ and 
$\Phi_{\rho,\tau,\gamma}(\omega_{\rho,\tau,\gamma}(x))=x.$
This readily implies that $\omega_{\rho,\tau,\gamma}(\epsilon)>\omega_{\rho,\tau,\gamma}(0)$ and therefore, since $g_\tau$ is strictly decreasing on $]-\infty; \omega_{\rho,\tau,\gamma}(\epsilon)]$, that $$g_\tau( \omega_{\rho,\tau,\gamma}(0)) >g_\tau( \omega_{\rho,\tau,\gamma}( \epsilon) )>-\frac{1}{\rho^2\gamma}.$$
Moreover,  by continuity, we readily have  \begin{equation}\label{der} \forall x \in ]-\epsilon;\epsilon[,~~\Phi_{\rho,\tau,\gamma}(\omega_{\rho,\tau,\gamma}(x))=x\end{equation}  and 
$ \omega_{\rho,\tau,\gamma}^{'}(0)\geq 0$ and $\Phi_{\rho,\tau,\gamma}^{'}( \omega_{\rho,\tau,\gamma}(0)) \geq 0.$
Now, differentiating both sides of (\ref{der}) at zero yields 
$\Phi_{\rho,\tau,\gamma}^{'}( \omega_{\rho,\tau,\gamma}(0)) \neq 0$ and $\omega_{\rho,\tau,\gamma}^{'}(0)\neq 0$ and then $\Phi_{\rho,\tau,\gamma}^{'}( \omega_{\rho,\tau,\gamma}(0)) > 0$ and  $\omega_{\rho,\tau,\gamma}^{'}(0)>0$. 

\noindent Thus, we have $\forall x \in   \mathbb{R}\setminus  \mbox{supp}(\mu_{\rho,\tau,\gamma})$, 
$ \omega_{\rho,\tau,\gamma}(x) \in  \mathbb{R}\setminus \mbox{supp}(\tau), \omega_{\rho,\tau,\gamma}^{'}(x)>~0, \\  ~~\Phi_{\rho,\tau,\gamma}^{'}( \omega_{\rho,\tau,\gamma}(x))>0, $~~$g_\tau( \omega_{\rho,\tau,\gamma}(x)) >-\frac{1}{\rho^2\gamma}$ and 
$$\Phi_{\rho,\tau,\gamma}(\omega_{\rho,\tau,\gamma}(x))=x.$$
These last  results combined with (ii) of Proposition \ref{inverserec} lead to  Proposition \ref{supportct}.\end{proof}
\subsection{Strict monotonicity of $\Phi_{\rho,\tau,\gamma}$  on ${\cal E}_{\rho,\tau,\gamma}$}\label{1}

\subsubsection{Strict monotonicity of $\Phi_{\rho,\tau,1}$  on ${\cal E}_{\rho,\tau,1}$}

The following proposition points out that $\Phi_{\rho,\tau,1}$ is globally increasing on ${\cal E}_{\rho,\tau,1}$.
\begin{proposition} \label{phicroitglobalement1}
For any $x_2 > x_1  $ in ${\cal E}_{\rho,\tau,1}$,
\begin{equation}\label{global}\Phi_{\rho,\tau,1}(x_2) > \Phi_{\rho,\tau,1}(x_1).\end{equation}
\end{proposition}
\begin{proof} Note that the result is obvious if $x_1=0$ or $x_2 =0$.
Let $\alpha$ be the symmetrization of the pushforward of $\tau$ by the map $t\mapsto \sqrt{t}$.
  Define for any $x\in \mathbb{R}\setminus \rm{supp}(\alpha)$,
\begin{eqnarray*}
H_{\rho ,\tau}(x) &=&x + \rho^2  g_\alpha(x)\\
&=& x + \rho^2  x g_\tau(x^2).
\end{eqnarray*}
 Let $x$ be in ${\cal E}_{\rho,\tau,1}$ such that $x>0$.
Note that $\sqrt{x}$ is in $\mathbb{R}\setminus \rm{supp}(\alpha)$ and 
\begin{equation}\label{phih}\Phi_{\rho,\tau,1}(x)=( H_{\rho,\tau}(\sqrt{x}))^2, \end{equation}
$$\Phi^{'}_{\rho,\tau,1}(x)= (1+\rho^2 g_\tau(x))H^{'}_{\rho ,\tau} (\sqrt{x}).$$
It readily follows that \begin{equation}\label{pos}1- \rho^2 \int \frac{d\alpha(t)}{(\sqrt{x}-t)^2}    =    H^{'}_{\rho ,\tau} (\sqrt{x})>0. \end{equation} 
Firstly, let $x_2 > x_1 > 0$ be in ${\cal E}_{\rho,\tau,1}$. Then  for $i=1,2$, $\sqrt{x_i} \in \mathbb{R}\setminus \rm{supp}(\alpha)$ and by (\ref{phih}) and (\ref{pos}), $\Phi_{\rho ,\tau,1}(x_i)=( H_{\rho ,\tau}(\sqrt{x_i}))^2,$
and
\begin{eqnarray*}
\frac{H_{\rho ,\tau }(\sqrt{x_2})-H_{\rho ,\tau }(\sqrt{x_1}) }{\sqrt{x_2}-\sqrt{x_1}}& = & 1 - \rho^2 \int_{\R}\frac{d\alpha (t)}{(\sqrt{x_1}-t)(\sqrt{x_2}-t)}  \\
& \geq &  1 -\left \{ \rho^2 \int \frac{d\alpha(t)}{(\sqrt{x_1}-t)^2}\right\}^{\frac{1}{2}}\left \{ \rho^2 \int \frac{d\alpha(t)}{(\sqrt{x_2}-t)^2}\right\}^{\frac{1}{2}}\\
&>&0.
\end{eqnarray*}
Since moreover  for $i=1,2$, $H_{\rho ,\tau}(\sqrt{x_i})=\sqrt{x_i} (1+ \rho^2  g_\tau(x_i))$ is positive, it follows that $$\Phi_{\rho ,\tau,1}(x_2)>\Phi_{\rho ,\tau,1}(x_1).$$
Secondly, let $ x_1<x_2<0 $ be in ${\cal E}_{\rho,\tau,1}$. Since $1+ \rho^2  g_\tau$ is strictly decreasing on $]-\infty; \min\rm{supp}(\tau) [$, we have $\forall x \in ]-\infty; x_2], ~ 1+ \rho^2  g_\tau(x)>0$ and  $\Phi^{'}_{\rho,\tau,1}(x)=  (1+\rho^2 g_\tau(x)) (2x  \rho^2 g_\tau^{'}(x) +  (1+\rho^2 g_\tau(x)))$ is obviously 
positive,  it follows that $$\Phi_{\rho,\tau,1}(x_2)>\Phi_{\rho,\tau,1}(x_1).$$
Finally, let $ x_1< 0<x_2 $ be in ${\cal E}_{\rho,\tau,1}$. We have   $\Phi_{\rho,\tau,1}(x_2)>0$ and  $\Phi_{\rho,\tau,1}(x_1)< 0$ so that 
$$\Phi_{\rho,\tau,1}(x_2)>\Phi_{\rho,\tau,1}(x_1).$$
Therefore, \begin{equation}\label{croitphi} \forall (x_1,x_2) \in \left({\cal E}_{\rho,\tau,1}\right)^2, \mbox{~if~} x_1<x_2, \mbox{~then~} \Phi_{\rho,\tau,1}(x_2)>\Phi_{\rho,\tau,1}(x_1).\end{equation} \end{proof}

\subsubsection{A relationship between  the  complements of the supports of $\mu_{\rho,\tau,\gamma}$
and    $\mu_{\rho \sqrt{\gamma},\tau,1}$}
The following lemmas will allow to establish,  in Proposition \ref{passagede1ac}, a relationship between  the  complement of the support of $\mu_{\rho,\tau,\gamma}$
and   the  complement of the support of $\mu_{\rho \sqrt{\gamma},\tau,1}$.

\begin{lemme}\label{incl1}
Recall that ${\cal E}_{\rho,\tau,\gamma}=\left\{ x \in  \mathbb{R}\setminus \mbox{supp}(\tau),  \Phi_{\rho,\tau,\gamma}^{'}(x) >0, g_\tau(x) >-\frac{1}{\rho^2 \gamma}\right\}$. We have for any $0<\gamma<1$, ${\cal E}_{\rho,\tau,\gamma} \subset {\cal E}_{\rho \sqrt{\gamma},\tau,1}$.
\end{lemme}
\noindent \begin{proof} This readily follows from the fact that $\Phi_{\rho,\tau,\gamma}^{'}(x)=\Phi_{\rho\sqrt{\gamma},\tau,1}^{'}(x)+ \gamma\rho^4 (1-\gamma) g^{'}_\tau (x)$ and then $\Phi_{\rho,\tau,\gamma}^{'}(x)<\Phi_{\rho \sqrt{\gamma},\tau,1}^{'}(x). $ \end{proof}
\begin{lemme} The function $K$ defined by 
$$K(x)=x +\frac{\rho^2(1-\gamma)}{1-\rho^2 \gamma g_{\mu_{\rho \sqrt{\gamma},\tau,1}}(x)}$$
is well defined   on $ \R\setminus \rm{supp}(\mu_{\rho \sqrt{\gamma},\tau,1})$ and 
\begin{equation}\label{comp} \forall x \in {\cal E}_{\rho \sqrt{\gamma},\tau,1}, ~~   \Phi_{\rho,\tau,\gamma}(x)=K(\Phi_{\rho \sqrt{\gamma},\tau,1}(x)).\end{equation}
\end{lemme}
\begin{proof}
 According to (\ref{Re}),
for any $x\neq 0$ in $\R\setminus \rm{supp}(\mu_{\rho \sqrt{\gamma},\tau,1})$, 
$\frac{1}{\rho^2 \gamma}-g_{\mu_{\rho \sqrt{\gamma},\tau,1}}(x) >0.$ Now, if $0\in \R\setminus \rm{supp}(\mu_{\rho \sqrt{\gamma},\tau,1})$, since $\frac{1}{\rho^2 \gamma}-g_{\mu_{\rho \sqrt{\gamma},\tau,1}}$ is strictly increasing on $]-\infty; \min  \rm{supp}(\mu_{\rho \sqrt{\gamma},\tau,1})[$, we can deduce that 
$\frac{1}{\rho^2 \gamma}-g_{\mu_{\rho \sqrt{\gamma},\tau,1}}(0) >0.$ Therefore for any $x$ in $\R\setminus \rm{supp}(\mu_{\rho \sqrt{\gamma},\tau,1})$, $\frac{1}{\rho^2 \gamma}-g_{\mu_{\rho \sqrt{\gamma},\tau,1}}(x) >0$ and $K$ is well defined.
Now, (\ref{comp}) readily follows from (\ref{TS3bis}) and (\ref{omegacircphi}) for $x\neq 0$ and may be proved for $x=0$ whenever $0 \in {\cal E}_{\rho \sqrt{\gamma},\tau,1}$ by a continuity argument. \end{proof}
\begin{proposition}\label{passagede1ac}
 $$\mathbb{R}\setminus  \mbox{supp}(\mu_{\rho,\tau,\gamma})=K\left(\left\{x \in \R\setminus \mbox{supp}(\mu_{\rho \sqrt{\gamma},\tau,1}), K'(x)>0 \right\} \right).$$
\end{proposition}
\begin{proof} Let $x$ be in $\mathbb{R}\setminus  \mbox{supp}(\mu_{\rho,\tau,\gamma})$. According to Proposition \ref{supportct}, there exists $u \in {\cal E}_{\rho,\tau,\gamma}$ such that $x=\Phi_{\rho,\tau,\gamma}(u)$. According to Lemma \ref{incl1}, $u$ belongs to $ {\cal E}_{\rho \sqrt{\gamma},\tau,1}$
and according to (\ref{comp}), $x=K(\Phi_{\rho \sqrt{\gamma},\tau,1}(u)).$ Proposition \ref{inverserec} (ii) implies that $\Phi_{\rho \sqrt{\gamma},\tau,1}(u)$ belongs to $\R\setminus \mbox{supp}(\mu_{\rho \sqrt{\gamma},\tau,1})$. Differentiating both sides of (\ref{comp}) leads to 
\begin{equation}\label{compder} \forall v \in {\cal E}_{\rho \sqrt{\gamma},\tau,1}, ~~   \Phi_{\rho,\tau,\gamma}^{'}(v)=K^{'}(\Phi_{\rho\sqrt{\gamma},\tau,1}(v))\Phi^{'}_{\rho\sqrt{\gamma},\tau,1}(v).\end{equation}
Then $K^{'}(\Phi_{\rho\sqrt{\gamma},\tau,1}(u))>0$ since $u \in {\cal E}_{\rho,\tau,\gamma} \subset {\cal E}_{\rho \sqrt{\gamma},\tau,1}$.
Therefore $$\mathbb{R}\setminus  \mbox{supp}(\mu_{\rho,\tau,\gamma})\subset K\left(\left\{x \in \R\setminus \mbox{supp}(\mu_{\rho \sqrt{\gamma},\tau,1}), K'(x)>0 \right\} \right).$$
Now, let $x$ be in $\R\setminus \mbox{supp}(\mu_{\rho \sqrt{\gamma},\tau,1})$ such that $K^{'}(x)>0$. According to Proposition \ref{supportct}, there exists $u$ in $ {\cal E}_{\rho \sqrt{\gamma},\tau,1}$ such that $x=\Phi_{\rho \sqrt{\gamma},\tau,1}(u)$
and (\ref{comp}) implies that $K(x)= K( \Phi_{\rho \sqrt{\gamma},\tau,1}(u))=\Phi_{\rho,\tau,\gamma}(u).$ Now, (\ref{compder}) implies that $\Phi^{'}_{\rho,\tau,\gamma}(u)>0$ and then $u\in {\cal E}_{\rho,\tau,\gamma}$ and, using Proposition \ref{inverserec} (ii), 
$\Phi_{\rho,\tau,\gamma}(u) \in \mathbb{R}\setminus  \mbox{supp}(\mu_{\rho,\tau,\gamma}).$  Therefore $K\left(\left\{x \in \R\setminus \mbox{supp}(\mu_{\rho \sqrt{\gamma},\tau,1}), K'(x)>0 \right\} \right) \subset \mathbb{R}\setminus  \mbox{supp}(\mu_{\rho,\tau,\gamma})$ and the proof is complete. \end{proof}
\subsubsection{Strict monotonicity of $\Phi_{\rho,\tau,\gamma}$  on ${\cal E}_{\rho,\tau,\gamma}$, when $\gamma<1$}

The following lemma will allow us to establish in Proposition \ref{phicroitglobalement} that $\Phi_{\rho,\tau,\gamma}$
is  globally increasing on the set $\left\{ x \in  \mathbb{R}\setminus \mbox{supp}(\tau), \Phi_{\rho,\tau,\gamma}^{'}(x) >0, g_\tau(x) >-\frac{1}{\rho^2\gamma}\right\}$.
\begin{lemme}\label{croitK}
Let $x_1<x_2$ be  in $\R\setminus  \rm{supp}(\mu_{\rho \sqrt{\gamma},\tau,1}) $ such that  for $ i=1,2$, $K^{'}(x_i)>0$. Then $K(x_1)<K(x_2)$.
\end{lemme}
\begin{proof} We have  for $ i=1,2$, $$K^{'}(x_i)= 1- \frac{\rho^4(1-\gamma)  \gamma}{(1-\rho^2 \gamma g_{\mu_{\rho \sqrt{\gamma},\tau,1}}(x_i))^2} \int \frac{d\mu_{\rho \sqrt{\gamma},\tau,1}(t)}{(x_i-t)^2}>0.$$
We have\\

\noindent $\frac{K(x_2)-K(x_1)}{ (x_2-x_1)}$ $$ = 1- \frac{\rho^4(1-\gamma)  \gamma}{(1-\rho^2 \gamma g_{\mu_{\rho \sqrt{\gamma},\tau,1}}(x_2))(1-\rho^2 \gamma g_{\mu_{\rho \sqrt{\gamma},\tau,1}}(x_1))} \int  \frac{d\mu_{\rho \sqrt{\gamma},\tau,1}(t)}{(x_2-t)(x_1-t)}.$$

\noindent Using Cauchy Schwartz inequality, we have \\

\noindent $\frac{\rho^4(1-\gamma)  \gamma}{(1-\rho^2 \gamma g_{\mu_{\rho \sqrt{\gamma},\tau,1}}(x_2))(1-\rho^2 \gamma g_{\mu_{\rho \sqrt{\gamma},\tau,1}}(x_1))} \int  \frac{d\mu_{\rho \sqrt{\gamma},\tau,1}(t)}{(x_2-t)(x_1-t)}
$
\begin{eqnarray*}&\leq & \left\{\frac{\rho^4(1-\gamma)  \gamma}{(1-\rho^2 \gamma g_{\mu_{\rho \sqrt{\gamma},\tau,1}}(x_1))^2} \int \frac{d\mu_{\rho \sqrt{\gamma},\tau,1}(t)}{(x_1-t)^2}\right\}^{\frac{1}{2}}\\
&& ~~~\times
\left\{\frac{\rho^4(1-\gamma)  \gamma}{(1-\rho^2  \gamma g_{\mu_{\rho \sqrt{\gamma},\tau,1}}(x_2))^2} \int \frac{d\mu_{\rho \sqrt{\gamma},\tau,1}(t)}{(x_2-t)^2}\right\}^{\frac{1}{2}}\\ &<&1.
\end{eqnarray*}
It follows that $K(x_2)-K(x_1)>0.$ \end{proof}

\begin{proposition}\label{phicroitglobalement}
For any $x_2 > x_1  $ in ${\cal E}_{\rho,\tau,\gamma}$,
\begin{equation}\label{global}\Phi_{\rho,\tau,\gamma}(x_2) > \Phi_{\rho,\tau,\gamma}(x_1).\end{equation}
\end{proposition}
\begin{proof} 
(\ref{compder}) implies that for any $x \in {\cal E}_{\rho,\tau,\gamma} \subset {\cal E}_{\rho \sqrt{\gamma},\tau,1}$, $K^{'} \left( \Phi_{\rho\sqrt{\gamma},\tau,1}(x)\right) >0$.
The proposition   readily follows from (\ref{comp}),  Lemma \ref{croitK} and Proposition \ref{phicroitglobalement1}.
 \end{proof}

\subsection{Proof of Theorem \ref{caractfinale} }\label{conc}



Propositions \ref{zeroinclus1} and  \ref{zeroinclus} correspond to A)  of Theorem \ref{caractfinale}.
B)  of Theorem \ref{caractfinale} can be easily deduced from 
(\ref{omegacircphi}), Propositions   \ref{supportct},
  \ref{phicroitglobalement1} and \ref{phicroitglobalement}. 
In order to establish  C) and D) of Theorem \ref{caractfinale}, we need to prove  preliminary lemmas.

\begin{lemme}\label{dd}
Let $I=]a;b[ $ be a bounded open interval and $F: I \rightarrow \R$ be a real analytic function   satisfying the following property ($\mathcal{F}$):
if $x<y $ are elements of $I$ such that $F^{'}(x)>0$ and  $F^{'}(y)>0$  then $F(x) < F(y)$.
 If $x_1$ and $y_1$ in  I satisfy, $ x_1<y_1, F^{'}(x_1)>0, F^{'}(y_1)<0$, then for all $x>y_1$, $x\in I$, we have $F^{'}(x) \leq 0.$ \end{lemme}
\begin{proof}
Set $ x_2= \min\{x>x_1, x\in I, F^{'}(x)<0\}$.     We have  $F^{'}(x_2)=0$. Since $F^{'}$ is analytic on $ I$, its zeroes are isolated. Therefore,
there exists $\epsilon>0$ small enough such that  $F^{'} <0$ on $ ]x_2,x_2+\epsilon]$.
Thus, $F(x_2+\epsilon)< F(x_2)$. Moreover, we have $x_2>x_1$ and $F^{'} \geq 0$ on $]x_1;x_2[$. Using once more the isolated zeroes principle, we can find $\epsilon^{'}$ small enough such that 
$F(x_2+\epsilon)<F(x_2-\epsilon^{'})<F(x_2)$  and $ F^{'}(x_2-\epsilon^{'})>0.$
Assume that there exists $ x >x_2+\epsilon$ in $ I $ such that $ F^{'}(x)>0$. Then, set $x_3= \min \{ x>x_2+\epsilon, x\in I, F^{'}(x)>0\}$.
We have $ F^{'}(x_3)=0$ and $F^{'}\leq 0$ on $ ]x_2+\epsilon, x_3[$. Thus  $F(x_3)\leq F(x_2+\epsilon)<F(x_2-\epsilon^{'})$. Using similar arguments as above 
we can find $\epsilon^{''}$  small enough such that $ F^{'}(x_3+\epsilon^{''})>0$ and $ F(x_3)<F(x_3+\epsilon^{''})<F(x_2-\epsilon^{'})$ which leads to a contradiction with property ($\mathcal{F}$).
Therefore there does not exist any $ x >x_2+\epsilon$ in $ I$ such that  $F^{'}(x)>0$. Moreover $F^{'} <0 $ on $ ]x_2, x_2+\epsilon] $. Hence we have $ F^{'}(x)\leq 0 $ for any $x>x_2$ and Lemma \ref{dd} follows. \end{proof}

The previous  lemma will be useful to establish the following one.
\begin{lemme}\label{dd2}
Let us consider an   open interval  I which is included in the set $\left\{x \in \mathbb{R} \setminus \mbox{supp}(\tau),   g_\tau(x) >-\frac{1}{\rho^2 \gamma}\right\}.$
If the set $\left\{x \in I,  \Phi_{\rho,\tau,\gamma}^{'}(x) >0\right\}$ is not empty then it is connected.
\end{lemme}
\begin{proof}   Assume first that $I$ is bounded: $I=]a;b[$ for some $a$ and $b$ in $\mathbb{R}$.
Set ${\cal O}=\{x \in I, \Phi_{\rho,\tau,\gamma}^{'}(x)>0\}$.
Assume that ${\cal O}$ is nonempty.
Define $$x_1=\min  \{x \in I, \Phi_{\rho,\tau,\gamma}^{'}(x)>0\},$$
$$x_2=\left\{\begin{array} {ll}\min\{x \in I, x>x_1, \Phi_{\rho,\tau,\gamma}^{'}(x)<0\}~\mbox{if~}   \{x \in I, x>x_1, \Phi_{\rho,\tau,\gamma}^{'}(x)<0\}\neq \emptyset,\\b \mbox{~else}.\end{array} \right.$$
Assume that $x_1=x_2$ which implies $x_2\neq b$ and then  $ \{x \in I, x>x_1, \Phi_{\rho,\tau,\gamma}^{'}(x)<0\}\neq \emptyset$ . For any   $\epsilon>0$ small enough, there exists $t_3$ in  $]x_1;  x_1+\epsilon[$ such that  $\Phi_{\rho,\tau,\gamma}^{'}(t_3)>0$.
Now, there exists 
$t_2$ in $ ]x_1,t_3[$, such that  $ \Phi_{\rho,\tau,\gamma}^{'}(t_2)<0$
and there exists $t_1$ in $]x_1,t_2[$ such that $ \Phi_{\rho,\tau,\gamma}^{'}(t_1)>0$,
which leads to a contradiction with  Lemma \ref{dd} (using Propositions \ref{phicroitglobalement1} and  \ref{phicroitglobalement}). 
Hence we must have $x_1<x_2$
and  Lemma \ref{dd} readily yields that $\Phi_{\rho,\tau,\gamma}^{'}(x) \leq 0$ for any $x$ in $I$ such that $x>x_2$ if $x_2<b$.\\
Therefore ${\cal O}=]x_1,x_2[ \setminus \{ x  \in ]x_1;x_2[, \Phi_{\rho,\tau,\gamma}^{'}(x)=0\}$. Now, using the isolated zeroes principle, B) of  Theorem \ref{caractfinale} and Propositions \ref{phicroitglobalement} and \ref{phicroitglobalement1}, it is easy to see that  a zero of $\Phi_{\rho,\tau,\gamma}^{'}$ in $ ]x_1;x_2[$ would yield a  mass for $\mu_{\rho,\tau,\gamma}$; therefore 
${\cal O}=]x_1;x_2[$.\\
Since $\lim_{x\rightarrow \pm \infty} \Phi_{\rho,\tau,\gamma}^{'}(x)= 1$ and  $\lim_{x\rightarrow \pm \infty} g_\tau(x)= 0$  , there exists $R>0$ such that $R> \max  \mbox{supp}(\tau)$  and  $g_\tau >-\frac{1}{\rho^2 \gamma}$,  $\Phi_{\rho,\tau,\gamma}^{'}>0$ on $]-\infty;-R]$ and $[R; +\infty[$.
Now, if $I$ is a  connected component  of the form $]-\infty; a[$ or $]a;+\infty[$, the previous  study can be carried out for $]-R;a[$ or $]a;R[$
and it readily follows that  ${\{x\in I , \Phi_{\rho,\tau,\gamma}^{'}(x)>0\}}$ is connected. 
 $\Box$\\

C) of Theorem \ref{caractfinale} corresponds to the following  lemma.

\begin{lemme}\label{CSomega} Let $a$ be in $\mathbb{R} \setminus \mbox{supp}(\mu_{\rho,\tau,\gamma}) $ such that $\omega_{\rho,\tau,\gamma}(a) \leq 0.$ Then
$a$ is to the  left  of $\mbox{supp}(\mu_{\rho,\tau,\gamma})$.
\end{lemme}
\noindent {\bf Proof:} First note that according to B)  of  Theorem \ref{caractfinale}, $\omega_{\rho,\tau,\gamma}(a) \in {\cal E}_{\rho,\tau,\gamma}$ and therefore 
$g_\tau(\omega_{\rho,\tau,\gamma}(a))>-\frac{1}{\rho^2\gamma}$ and $ \Phi_{\rho,\tau,\gamma}^{'}(\omega_{\rho,\tau,\gamma}(a))>0$. Since $g_\tau$ is stricly deacreasing on $]-\infty; \min \mbox{supp}(\tau)[$, we can deduce that $$]-\infty; \omega_{\rho,\tau,\gamma}(a)] \subset \left\{ x \in  \mathbb{R}\setminus \mbox{supp}(\tau),   g_\tau(x) >-\frac{1}{\rho^2\gamma}\right\}.$$
Now, since $\lim_{x\rightarrow -\infty} \Phi_{\rho,\tau,\gamma}^{'}(x)= 1$ and $\Phi_{\rho,\tau,\gamma}^{'}(\omega_{\rho,\tau,\gamma}(a)) >0$, Lemma \ref{dd2}
readily implies that $]-\infty; \omega_{\rho,\tau,\gamma}(a)] \subset  {\cal E}_{\rho,\tau,\gamma}$ and then, using B)  of Theorem  \ref{caractfinale},
$]-\infty;a] \subset \mathbb{R} \setminus \mbox{supp}(\mu_{\rho,\tau,\gamma}). $\end{proof}


When  the support of $\tau$ has a finite number of connected components, we have the following description of the support of $\mu_{\rho,\tau,\gamma}$ in terms of a finite union of closed disjoint intervals.
\begin{proposition}\label{compsup}
For any $0<\gamma\leq 1$, recall that  $$ {\cal E}_{\rho,\tau,\gamma}=\left\{ x \in  \mathbb{R}\setminus \mbox{supp}(\tau),  \Phi_{\rho,\tau,\gamma}^{'}(x) >0, g_\tau(x) >-\frac{1}{\rho^2\gamma}\right\}.$$ 
Assume that the support of $\tau$ is a finite union of disjoint (possibly degenerate) closed bounded intervals. There exists a nonnul integer number  $p$ and  $u_1< v_1<u_2<\ldots <u_p < v_p$ depending on $\rho,\tau,\gamma $, such that \begin{equation}\label{vrai} {{\cal E}_{\rho,\tau,\gamma}}=]-\infty;u_1[ \cup_{l=1}^{p-1}]v_l;u_{l+1}[\cup 
]v_p;+\infty[. \end{equation}
We have \begin{equation}\label{vrai2}\mbox{supp}(\tau) \subset \cup_{l=1}^{p} [u_l;v_l]\end{equation}
and for each $l \in \{1,\ldots,p\}$,   \begin{equation}\label{vrai3}[u_l;v_l]\cap \mbox{supp}(\tau) \neq \emptyset.\end{equation} 
Moreover,
\begin{equation}\label{compcon}
 \mbox{supp}(\mu_{\rho,\tau,\gamma})
=\cup_{l=1}^p [\Phi_{\rho,\tau,\gamma}(u_l^-);\Phi_{\rho,\tau,\gamma}(v_l^+)],\end{equation}
with
 for all $ l=1,\ldots, p-1,  \Phi_{\rho,\tau,\gamma}(v_l^+) <  \Phi_{\rho,\tau,\gamma}(u_{l+1}^-),$ and 
 for all $ l=1,\ldots,p, $ \begin{equation}\label{ext} \Phi_{\rho,\tau,\gamma}(u_l^-) < \Phi_{\rho,\tau,\gamma}(v_l^+).\end{equation}
\end{proposition}
\begin{proof}
Since 
$1+ \rho^2 \gamma g_\tau$ is strictly decreasing on each connected component of $ \mathbb{R}\setminus \mbox{supp}(\tau)$ and $\lim_{x \rightarrow \pm \infty}
g_\tau(x)=0$, it is easy to see that 
there exists a nonnul integer number  $m$ and  $a_1< b_1<a_2<\ldots <a_m< b_m$ depending on $\rho,\tau,\gamma $, such that \begin{equation}\label{vrai1} \left\{ x \in  \mathbb{R}\setminus \mbox{supp}(\tau),   g_\tau(x) >-\frac{1}{\rho^2 \gamma}\right\}=]-\infty;a_1[ \cup_{l=1}^{m-1}]b_l;a_{l+1}[\cup 
]b_m;+\infty[, \end{equation}
$$\mbox{supp}(\tau) \subset \cup_{l=1}^{m} [a_l;b_l]$$
and for each $l \in \{1,\ldots,m\}$,   $[a_l;b_l]\cap \mbox{supp}(\tau) \neq \emptyset$ .
Lemma \ref{dd2} (using also  that  $\lim_{x\rightarrow \pm \infty} \Phi_{\rho,\tau,\gamma}^{'}(x)= 1$) readily implies (\ref{vrai}), (\ref{vrai2}) and (\ref{vrai3}).
(\ref{compcon}) can be easily deduced from  B)   of Theorem \ref{caractfinale}. Since $\mu_{\rho,\tau,\gamma}$ has no mass, (\ref{ext}) follows.
The proof of Proposition \ref{compsup} is complete.
 \end{proof}



Proposition \ref{compsup} corresponds to D) of Theorem \ref{caractfinale}. Note that choosing a matricial model as introduced in Section \ref{spi} but without spikes, it is easy to deduce (\ref{palier}) from 
Theorem \ref{sep}, Lemma \ref{inclusionN},  Theorem \ref{inclusion} A) and the weak convergence almost surely  of the spectral measures. The proof of Theorem \ref{caractfinale} is complete.

\section{Exact separation phenomenon}\label{thsep}
\subsection{Preliminary results}
Our proof of Theorem \ref{sep} will need  to establish the following preliminary  proposition and lemmas.
\begin{proposition}\label{inc}:
Let $\gamma $ be in $]0;1]$, $\rho$ be in $]0; +\infty[$ and $\tau$ be a compactly supported probability measure on $[0;+\infty[$. 
For all $0< \hat \rho < \rho $, ${\cal E}_{\rho,\tau,\gamma }\subset \, 
{\cal E}_{\hat \rho,\tau,\gamma  }$ 
so that it makes sense to consider the following composition of homeomorphisms 
$$\Phi_{\hat \rho , \tau ,\gamma}\circ \omega_{\rho,\tau,\gamma}: \, 
   \mathbb{R}\setminus \rm{supp}(\mu_{ \rho,\tau,\gamma})   \rightarrow 
\Phi_{\hat \rho , \tau ,\gamma}({\cal E}_{\rho , \tau,\gamma } )\subset \, 
  \mathbb{R}\setminus  \rm{supp}(\mu_{  \hat{\rho},\tau,\gamma}) ,$$
which is stricly increasing on 
$  \mathbb{R}\setminus  \rm{supp}(\mu_{ \rho,\tau,\gamma})  $.
\end{proposition}
\begin{proof} Let  $0< \hat \rho < \rho $. Let $u$ be in ${\cal E}_{\rho,\tau,\gamma}$. $g_\tau(u) > -\frac{1}{\rho^2 \gamma}$ obviously implies $g_\tau(u) > -\frac{1}{\hat{\rho}^2 \gamma}$.
We have  $$\Phi^{'}_{ \hat{\rho},\tau,\gamma} (u)= (1+\gamma \hat{\rho}^2g_{ \tau}(u))\left\{ (1+\gamma \hat{\rho}^2g_{ \tau}(u))+ \left[2u \gamma \hat{\rho}^2    +   \frac{\hat{\rho}^2 (1-\gamma) \gamma \hat{\rho}^2}{1+\gamma \hat{\rho}^2g_{ \tau}(u)}\right] g^{'}_\tau(u) \right\}.$$
\begin{itemize}
\item If $ 2u \gamma \hat{\rho}^2    +   \frac{\hat{\rho}^2 (1-\gamma) \gamma \hat{\rho}^2}{1+\gamma \hat{\rho}^2g_{ \tau}(u)}\leq 0$, it is obvious that $\Phi^{'}_{ \hat{\rho},\tau,\gamma} (u)>0$.
\item If $ 2u \gamma    +   \frac{ (1-\gamma) \gamma}{\frac{1}{\hat{\rho}^2}+\gamma g_{ \tau}(u)}>0$, then $ 2u \gamma    +   \frac{ (1-\gamma) \gamma }{\frac{1}{{\rho}^2}+\gamma g_{ \tau}(u)}>0$.
Moreover $\Phi^{'}_{\rho,\tau,\gamma}(u) >0$ means that 
$$ (1+\gamma \rho^2g_{ \tau}(u))^2 + 2u (1+\gamma \rho^2g_{ \tau}(u))\gamma \rho^2 g^{'}_{ \tau}(u)   +   \rho^2 (1-\gamma) \gamma  \rho^2 g^{'}_\tau(u) >0.$$
Since $(1+\gamma \rho^2g_{ \tau}(u))>0$ and  $ 2u \gamma    +   \frac{ (1-\gamma) \gamma }{\frac{1}{{\rho}^2}+\gamma g_{ \tau}(u)}>0$, this implies 
\begin{equation}\label{min}g^{'}_\tau (u) > - \frac{(1+\gamma \rho^2g_{ \tau}(u))}{2u\gamma \rho^2 + \frac{\rho^2(1-\gamma) \gamma \rho^2}{1+ \gamma \rho^2 g_\tau(u)}}. \end{equation}
(\ref{min}) implies that 
\begin{eqnarray*}&& (1+\gamma \hat{\rho}^2g_{ \tau}(u))+ \left[2u \gamma \hat{\rho}^2    +   \frac{\hat{\rho}^2 (1-\gamma) \gamma \hat{\rho}^2}{1+\gamma \hat{\rho}^2g_{ \tau}(u)}\right] g^{'}_\tau(u)\\
&>&  (1+\gamma \hat{\rho}^2g_{ \tau}(u)) -\frac{\hat{\rho}^2}{\rho^2}  \frac{(1+\gamma\rho^2g_{ \tau}(u))}{2u\gamma  + \frac{(1-\gamma) \gamma\rho^2}{1+ \gamma\rho^2 g_\tau(u)}} \left(2u\gamma  + \frac{(1-\gamma) \gamma\hat{\rho}^2}{1+ \gamma\hat{\rho}^2 g_\tau(u)}\right)\\
&>&  (1+\gamma\hat{\rho}^2g_{ \tau}(u))\\&& -\frac{\hat{\rho}^2}{\rho^2}  \frac{(1+\gamma \rho^2g_{ \tau}(u))}{2u\gamma   + \frac{(1-\gamma) \gamma\rho^2}{1+ \gamma \rho^2 g_\tau(u)}} \left(2u\gamma + \frac{(1-\gamma ) \gamma{\rho}^2}{1+ \gamma{\rho}^2 g_\tau(u)}
 \right. \\ && ~~~~~~~~~~~~~~~~~~~~~~~~~~~~~+\left. (1-\gamma) \gamma \left[\frac{\hat{\rho}^2}{1+ \gamma\hat{\rho}^2 g_\tau(u)} -\frac{{\rho}^2}{1+ \gamma{\rho}^2 g_\tau(u)}\right]\right)\\
&=& 1 -\frac{\hat{\rho}^2}{\rho^2}+ \Delta(u)\end{eqnarray*}
\noindent where \begin{eqnarray*} \Delta(u)&= & -(1-\gamma)\gamma\frac{\hat{\rho}^2}{\rho^2}  \frac{(1+\gamma\rho^2g_{ \tau}(u))}{2u\gamma  + \frac{(1-\gamma) \gamma\rho^2}{1+ c\rho^2 g_\tau(u)}}\left[\frac{\hat{\rho}^2}{1+ c\hat{\rho}^2 g_\tau(u)} -\frac{{\rho}^2}{1+ c{\rho}^2 g_\tau(u)}\right]\\
&=&(1-~\gamma)\gamma\frac{\hat{\rho}^2}{\rho^2}(\rho^2-\hat{\rho}^2)  \frac{(1+\gamma\rho^2g_{ \tau}(u))}{(1+\gamma\hat{\rho}^2g_{ \tau}(u))}\frac{1}{2u \gamma (1+\gamma\rho^2g_{ \tau}(u))+ (1-\gamma) \gamma \rho^2}\\&>&0
\end{eqnarray*}
It follows that $\Phi^{'}_{ \hat{\rho},\tau,\gamma}(u) >0$. 
\end{itemize}
$\Phi_{ \hat{\rho},\tau,\gamma}\circ \omega_{\rho,\tau,\gamma}$ is strictly increasing on  $ \mathbb{R}\setminus \rm{supp}(\mu_{ \rho,\tau,\gamma})  $ by B)  of Theorem \ref{caractfinale}. The proof is complete. \end{proof}

\begin{lemme}\label{convg} Let $[a,b]$ be such that there exists $\delta>0$ such that for all large $N$, $[a-\delta;b+\delta] \subset \R \setminus \rm{supp}( \mu_{\sigma, \mu_{A_NA_N^*},c_N})$.
Then for all $ 0< \tau < \delta$, for all large $N$, $$\omega_{\sigma, \nu, c}\left( [a-\tau; b+\tau] \right) =[\omega_{\sigma, \nu, c}(a-\tau); \omega_{\sigma, \nu, c}(b+\tau)] \subset \R \setminus \rm{supp}(\mu_{A_NA_N^*}).$$ Moreover, for any $x \in [a;b]$, when $N$ goes to infinity,
$g_{\mu_{A_NA_N^*}} \left( \omega_{\sigma, \mu_{A_NA_N^*}, c_N}(x)\right)$ converges towards $g_{\nu} \left( \omega_{\sigma, \nu, c}(x)\right)$.
\end{lemme}
\begin{proof} 
We will use  the obvious fact that if a sequence $\mu_n$ of probability measures weakly 
converges towards a probability measure $\mu$, if for some $\delta>0$, $[u-\delta; u+\delta] \subset \mathbb{R}\setminus \rm{supp}~{\mu_n}$ for all large $n$, then $g_{\mu_n}(u)$ converges towards $g_\mu(u)$.\\
Let $0< \tau< \tau^{'}<\delta$.
Since $ \mu_{\sigma, \mu_{A_NA_N^*},c_N}$ weakly converges towards $ \mu_{\sigma, \nu,c}$, for any $u \in [a-\tau^{'}; b+\tau^{'}]$, 
$ \omega_{\sigma, \mu_{A_NA_N^*}, c_N}(u)$ converges towards $ \omega_{\sigma, \nu, c}(u)$ and
$ [a-\tau^{'}; b+\tau^{'}] \subset \mathbb{R} \setminus \rm{supp} ( \mu_{\sigma, \nu, c})$. According to B)  of Theorem \ref{caractfinale} , $ \omega^{'}_{\sigma, \mu_{A_NA_N^*}, c_N}>0$ and $\omega^{'}_{\sigma, \nu, c}>0$ on $ [a-\tau^{'}; b+\tau^{'}]$  . It readily follows that, for all large $N$,
\begin{eqnarray*}[\omega_{\sigma, \nu, c}(a-\tau); \omega_{\sigma, \nu, c}(b+\tau)]& \subset& ] \omega_{\sigma, \mu_{A_NA_N^*}, c_N}(a-\tau^{'}); \omega_{\sigma, \mu_{A_NA_N^*}, c_N}(b+\tau^{'})[ \\&\subset& \R \setminus  \rm{supp}(\mu_{A_NA_N^*})\end{eqnarray*}
\noindent where we used B)  of Theorem \ref{caractfinale} for the last inclusion. Using Vitali's Theorem, we deduce that $g_{\mu_{A_N A_N^*}}$ converges uniformly on $[\omega_{\sigma, \nu, c}(a-\tau/2); \omega_{\sigma, \nu, c}(b+\tau/2)]$ and then that, for any $x \in [a;b]$,
%
%
 $ g_{\mu_{A_N A_N^*}}(\omega_{\sigma, \mu_{A_NA_N^*}, c_N}(x))$ converges towards $g_\nu(\omega_{\sigma, \nu, c}(x)).$ 
\end{proof}

\begin{lemme}\label{gap}
Let $[a,b]$  ($a\neq b$) be a compact set satisfying the couple ${\cal P}(\sigma)$ of properties   defined in Theorem \ref{pasde} 7. and such that, if $c<1$, $\omega_{\sigma,\nu,c}(a)> 0$. 
 For all $0 \leq  \hat \sigma \leq \sigma $, define $a(\hat \sigma)=
\Phi_{\hat \sigma , \nu,c }(\omega_{\sigma,\nu,c}(a))$ and $b(\hat \sigma)= \Phi_{\hat \sigma , \nu,c }(\omega_{\sigma,\nu,c}(b)).$
Then for each $0 <  \hat \sigma < \sigma $, $[a(\hat \sigma), b(\hat \sigma)]$ satisfies the couple ${\cal P}(\hat \sigma)$ of properties. 
Moreover, there exists $m_{a,b}>0$ such that for all $0 \leq \hat \sigma \leq \sigma $, 
\begin{equation}\label{minor} \sqrt{b(\hat \sigma)}- \sqrt{a(\hat \sigma)}\geq m_{a,b}.\end{equation}
\end{lemme}

\begin{proof} Denote by $\delta$ and $\tau$ the parameters introduced in 7. of Theorem \ref{pasde}.
Note that if $c=1$,  since $a>0$ and using (\ref{Re}), we have \begin{equation}\label{wpo} \omega_{\sigma , \nu,c }(a)>0. \end{equation}
Moreover, according to B)  of Theorem \ref{caractfinale} and Proposition \ref{inc}, for any $0<c\leq 1$ and any $0< \hat{\sigma} < \sigma$, we have  $1+ c\hat{\sigma}^2 g_\nu( \omega_{\sigma, \nu, c}(a))>0$.
Therefore  it readily follows that if $c=1$, then
$a(\hat \sigma) >0$ and if $c<1$, if $\omega_{\sigma,\nu,c}(a)> 0$, then $a(\hat \sigma) >0$.

Moreover, since according to Proposition \ref{inc}, $x \mapsto \Phi_{\hat \sigma , \nu ,c}\circ \omega_{\sigma,\nu,c}(x)$ is strictly increasing on $[a;b]$, we have \begin{equation}\label{ab} 0 < a(\hat{\sigma}) < b(\hat{\sigma}). \end{equation}
According to Lemma \ref{convg}, for all $0< \tau^{'}< \delta$ and for all  $ x \in[a-\tau^{'}; b+\tau^{'}]$, $g_{\mu_{A_NA_N^*}} \left( \omega_{\sigma, \mu_{A_NA_N^*}, c_N}(x)\right)$ converges towards $g_{\nu} \left( \omega_{\sigma, \nu, c}(x)\right)$. It readily follows that for any $\hat{\sigma} >0$,
 $\Phi_{\hat{\sigma}, \mu_{A_NA_N^*}, c_N} \left( \omega_{\sigma, \mu_{A_NA_N^*}, c_N}(x)\right)$ converges towards $\Phi_{\hat{\sigma}, \nu, c}  \left( \omega_{\sigma, \nu, c}(x)\right)$. According to Proposition  \ref{inc}, $x \mapsto \Phi_{\hat \sigma , \nu ,c}( \omega_{\sigma,\nu,c}(x))$ and 
$x \mapsto \Phi_{\hat \sigma ,  \mu_{A_NA_N^*}, c_N}( \omega_{\sigma, \mu_{A_NA_N^*}, c_N}(x))$ are strictly increasing on $[a-\tau^{'};b+\tau^{'}]$.
It readily follows that, for all $0< \tau^{''}<\tau^{'}< \delta$, for all large $N$,\\

$ [ \Phi_{\hat \sigma , \nu ,c}\circ \omega_{\sigma,\nu,c}(a-\tau^{''}); \Phi_{\hat \sigma , \nu ,c}\circ \omega_{\sigma,\nu,c}(b+\tau^{''})]$
\begin{eqnarray*}&\subset&  ] \Phi_{\hat \sigma ,  \mu_{A_NA_N^*}, c_N}( \omega_{\sigma, \mu_{A_NA_N^*}, c_N}(a-\tau^{'})); \Phi_{\hat \sigma ,  \mu_{A_NA_N^*}, c_N}( \omega_{\sigma, \mu_{A_NA_N^*}, c_N}(b+\tau^{'}))[\\
&\subset& \R \setminus \rm{supp}(\mu_{\hat \sigma ,  \mu_{A_NA_N^*}, c_N}).
\end{eqnarray*}
Therefore, there exists $\alpha >0$ such that $ \Phi_{\hat \sigma , \nu ,c}( \omega_{\sigma,\nu,c}(a))-\alpha >0$ and 
$$ [ \Phi_{\hat \sigma , \nu ,c}( \omega_{\sigma,\nu,c}(a))-\alpha; \Phi_{\hat \sigma , \nu ,c} (\omega_{\sigma,\nu,c}(b))+\alpha]\subset \R \setminus \rm{supp}(\mu_{\hat \sigma ,  \mu_{A_NA_N^*}, c_N})$$ \noindent which is the first property of ${\cal P}(\hat \sigma)$.\\

The second property of ${\cal P}(\hat \sigma)$ will be obviously satisfied if we prove that there exists $\epsilon >0$ such that for all large $N$,
\begin{equation}\label{etoile} 
\omega_{\hat \sigma ,  \mu_{A_NA_N^*}, c_N}(]a(\hat \sigma)-\epsilon; b(\hat \sigma)+\epsilon[ )\subset \omega_{ \sigma ,  \mu_{A_NA_N^*}, c_N}(]a-\tau; b+\tau[) .
\end{equation}
We have proved that there exists $\alpha>0$ such that \\

\noindent $ [ \Phi_{\hat \sigma , \nu ,c}( \omega_{\sigma,\nu,c}(a))-\alpha; \Phi_{\hat \sigma , \nu ,c}( \omega_{\sigma,\nu,c}(b))+\alpha]$
\begin{eqnarray*}&\subset&  ] \Phi_{\hat \sigma ,  \mu_{A_NA_N^*}, c_N}( \omega_{\sigma, \mu_{A_NA_N^*}, c_N}(a-\tau^{})); \Phi_{\hat \sigma ,  \mu_{A_NA_N^*}, c_N}( \omega_{\sigma, \mu_{A_NA_N^*}, c_N}(b+\tau^{}))[\\
&\subset& \R \setminus \rm{supp}(\mu_{\hat \sigma ,  \mu_{A_NA_N^*}, c_N}).
\end{eqnarray*}
(\ref{etoile}) follows from the fact  that, according to B) of Theorem \ref{caractfinale}, $ \omega_{\hat \sigma, \mu_{A_NA_N^*}, c_N}$ is strictly increasing on the intervals involved in the inclusion and moreover  $\omega_{\hat{\sigma},\mu_{A_NA_N^*},c_N}\circ \Phi_{\hat{\sigma}, \mu_{A_NA_N^*}, c_N}(x) =x$ for each $x$ in  
${\cal E}_{\hat{\sigma}, \mu_{A_NA_N^*}, c_N}$ .\\

\noindent  According to (\ref{ab}) (and B) of Theorem \ref{caractfinale} for  $\hat{\sigma}=\sigma$ or $ \hat{\sigma}=0$), we have $$\forall 0\leq \hat{\sigma} \leq \sigma, ~\sqrt{b(\hat{\sigma}) }-\sqrt{a(\hat{\sigma})}>0.$$ Since $\hat{\sigma} \mapsto \sqrt{b(\hat{\sigma})}-\sqrt{a(\hat{\sigma})}$ is obviously continuous on $[0;\sigma]$, 
 (\ref{minor}) readily follows. 
\end{proof}

The proof of Theorem \ref{sep} will also rely on the following classical result.

\begin{lemme}{(cf. Theorem 11.8 of \cite{BaiSilverbook})} \label{Weyl}
Let B and C be two $n \times N$ complex matrices, $n\leq N$. 
For any pair of integers $j,k$ such that $1\leq j,k\leq n$ and $j+k\leq n+1$, 
we have
$$\sqrt{\lambda_{j+k-1}[(B+C)(B+C)^*]}\leq \sqrt{\lambda_{j}(BB^*)}+\sqrt{\lambda_{k}(CC^*)}.$$
\end{lemme}
\subsection{ Proof of Theorem \ref{sep}}
First, we prove the following weaker proposition the proof of which is in the lineage of \cite{BaiSil99}, \cite{CDF09},\cite{CDFF}.
\begin{proposition}\label{seppre}Assume conditions [1-7] of Theorem \ref{pasde} are satisfied.
If $c <1$,  assume moreover that 
$\omega_{\sigma,\nu,c}(a)> 0$.
 Then for  $N$ large enough, $$\omega_{{\sigma,\nu,c}}([a,b])=[\omega_{{\sigma,\nu,c}}(a);\omega_{{\sigma,\nu,c}}(b)] \subset \R \setminus  \mbox{supp}(\mu _{A_N A_N^*}).$$
With the convention that $\lambda _0(M_N)=\lambda _0(A_NA_N^*)=+\infty $ and 
$\lambda _{N+1}(M_N)=\lambda _{N+1}(A_NA_N^*)=-\infty $, for $N$ large enough, let  $i_N\in \{0,\ldots,n\}$
be such that
\begin{equation}\label{iN2}\lambda_{i_N+1}(A_N A_N^*) <\omega_{{\sigma,\nu,c}}(a) \mbox{~~ and ~~} \lambda_{i_N}(A_N A_N^*) > \omega_{{\sigma,\nu ,c}}(b).\end{equation}
Then
$$P[\mbox{for all large N}, \lambda_{i_N+1}(M_N) <a\mbox{~and} ~ \lambda_{i_N}(M_N)> 
b] = 1.$$
\end{proposition}

\begin{proof} 
According to  (\ref{minor}) in Lemma \ref{gap}, there exists ${m}_{a,b}>0$ such that for all for all $\hat{\sigma} \in [0; \sigma]$,
\begin{equation}\label{largeur} \sqrt{ \Phi_{\hat{\sigma},\nu,c}(\omega_{\sigma , \nu,c }(b))}-\sqrt{ \Phi_{\hat{\sigma},\nu,c}(\omega_{\sigma , \nu,c }(a))}\geq {m}_{a,b}. \end{equation}
Since  $\hat{\sigma} \mapsto \sqrt{ \Phi_{\hat{\sigma},\nu,c} (\omega_{\sigma , \nu,c }(a))}$ and $\hat{\sigma} \mapsto \sqrt{\Phi_{\hat{\sigma},\nu,c} (\omega_{\sigma , \nu,c }(b))}$ are uniformly continuous  on $[0; \sigma]$, there exists \begin{equation} \label{delta}0< \delta ({a,b}) <\frac{{m}_{a,b}}{4} \end{equation}
such that for any $\hat{\sigma}$, $\hat{\sigma}^{'}$ in $[0; \sigma]$ satisfying $\vert \hat{\sigma}-\hat{\sigma}^{'}\vert \leq \delta({a,b})$, we have 
\begin{equation}\label{cont} \begin{array}{ll}\vert \sqrt{ \Phi_{\hat{\sigma},\nu,c}(\omega_{\sigma , \nu,c }(a))}-\sqrt{ \Phi_{\hat{\sigma}^{'},\nu,c}(\omega_{\sigma , \nu, c }(a))}\vert < \frac{{m}_{a,b}}{4},\\
\mbox{and ~} \vert \sqrt{ \Phi_{\hat{\sigma},\nu,c}(\omega_{\sigma , \nu,c }(b))}-\sqrt{ \Phi_{\hat{\sigma}^{'},\nu,c}(\omega_{\sigma , \nu,c }(b))}\vert < \frac{{m}_{a,b}}{4}.\end{array}\end{equation}
According to \cite{CSBD}, we can choose $C>1$ large enough such that almost surely, for all large $N$, 
\begin{equation}\label{max} 0 \leq \sqrt{\lambda _1(X_NX_N^*)} < C .\end{equation}
Then, let us choose \begin{equation}\label{Cab}0 < C_{a,b}< \frac{\delta({a,b})}{C\sigma}.\end{equation} Given $k\geq 0$ define 
$$\sigma _k=\sigma (1-\frac{1}{1+kC_{a,b}}),$$
 $$s_k={ \Phi_{{\sigma}_k,\nu,c}(\omega_{\sigma , \nu,c }(a))}$$
and  $$t_k={ \Phi_{{\sigma}_k,\nu,c}(\omega_{\sigma , \nu,c }(b))}.$$
Note that for all $k$,
\begin{equation}\label{k}\vert \sigma_{k+1} -\sigma_k\vert \leq C_{a,b} \sigma. \end{equation}
According to (\ref{largeur}), for any $k$, \begin{equation}\label{ouf}\sqrt{ t_k}-\sqrt{s_k} \geq {m}_{a,b}\end{equation}
and according to (\ref{cont}), (\ref{Cab}) and (\ref{k}), 
\begin{equation}\label{petitpas}\vert \sqrt{s_{k+1}}-\sqrt{s_k}\vert < \frac{{m}_{a,b}}{4} ~~\mbox{and ~}\vert \sqrt{t_{k+1}}-\sqrt{t_k}\vert < \frac{{m}_{a,b}}{4}.\end{equation}

\noindent Now, let us  introduce a sequence of matrices $M_N^{(k)}$ interpolating from $M_N$ to $A_NA_N^*$: 
$$M_N^{(k)}:=({\sigma _k}X_N+A_N)({\sigma _k}X_N+A_N)^*.$$
For all $k \geq 0$, 
set $${\rm E}_k=\{ \text{no eigenvalues of $M_N^{(k)}$ in $[s_{k},t_{k}]$, for all large $N$}\} .$$
By Lemma \ref{gap}  and Theorem \ref{pasde}, 
we know that $\mathbb P({\rm E}_k)=1$ for all $k>0$.  
Moreover, according to Lemma \ref{convg}, for $N\geq N_0$,  $[\omega_{\sigma, \nu, c}(a); \omega_{\sigma, \nu, c}(b)] \subset \R \setminus \rm{supp}(\mu_{A_NA_N^*}).$
For $N \geq N_0$, let
$i_N$  be such that 
\begin{equation}\label{iN}
\lambda _{i_N+1}(A_NA_N^*) < \omega_{\sigma, \nu, c}(a) \, \text{ and } \, 
\lambda _{i_N}(A_NA_N^*) > \omega_{\sigma, \nu, c}(b). \end{equation}
 We shall  show that by induction on  $k$ that one has for all $k \geq 0$, 
\begin{equation} {\label{casek}}
\mathbb P[\lambda _{i_N+1}(M_N^{(k)}) < s_k \, \text{ and } \, 
\lambda _{i_N}(M_N^{(k)}) > t_k, \, \text{ for all large $N$}]=1.
\end{equation}
This is true for $k=0$ since $\sigma_0=0$, $M_N^{(0)}=A_NA_N^*$ and $\Phi_{0,\nu,c}(x)=x.$ 
Now, let us assume that \eqref{casek} holds true for some $k \geq 0$. 
If $i_N=N$, it is obvious that $\lambda_{i_N+1}(M_N^{(k+1)}) < s_{k+1} $. If $i_N \in \{0;\ldots;N-1\}$,
we can deduce from Lemma \ref{Weyl} that 
$$\sqrt{\lambda _{i_N+1}(M_N^{(k+1)})}\leq \sqrt{\lambda _{i_N+1}(M_N^{(k)})}+\sqrt{\lambda _1(X_NX_N^*)}C_{a,b}\sigma.$$
It follows using \eqref{max}, \eqref{delta} and \eqref{Cab} that almost surely for all large $N$,
$$\sqrt{\lambda _{i_N+1}(M_N^{(k+1)})} < \sqrt{s_k} + \frac{{m}_{a,b}}{4}:=\hat s_k.$$
If $i_N=0$, $\lambda_{i_N}(M_N^{(k+1)}) > t_{k+1}$ and, if $i_N \in \{1;\ldots;N\}$,
similarly, one can show that  almost surely for all large $N$,
$$\sqrt{\lambda_{i_N}(M_N^{(k+1)})} > \sqrt{t_k}-\frac{{m}_{a,b}}{4}:=\hat t_k.$$
Inequalities \eqref{petitpas} and (\ref{ouf})  ensure that $$[\hat s_k, \hat t_k]\subset [\sqrt{s_{k+1}},\sqrt{ t_{k+1}}].$$
As $\mathbb P({\rm E}_{k+1})=1$, we deduce that, with probability $1$, 
$$\lambda_{i_N+1}(M_N^{(k+1)}) < s_{k+1} \, \text{ and } \, 
\lambda_{i_N}(M_N^{(k+1)}) > t_{k+1}, \quad \text{for all large $N$}.$$
This completes the proof by induction of \eqref{casek}.\\

\noindent Now, we are going to show that there exists $K$ large enough so that, 
for all $k \geq K$, there is an exact separation of the eigenvalues of the matrices $M_N$ and $M_N^{(k)}$.\\
Let $\alpha $ be such that $0< \alpha< \frac{\sqrt{b} -\sqrt{a}}{2}$ and 
  $0< \delta <\frac{\alpha}{2}$ be
such that for any $\hat{\sigma}$ in $[0; \sigma]$ satisfying $\vert \hat{\sigma}-{\sigma}\vert \leq \delta$, we have 
$$ \vert \sqrt{ \Phi_{\hat{\sigma},\nu,c}(\omega_{\sigma , \nu,c }(a))}-\sqrt{ a}\vert \leq \frac{\alpha}{2}
 \mbox{~~and~~}\vert \sqrt{ \Phi_{\hat{\sigma},\nu,c}(\omega_{\sigma , \nu,c }(b))}-\sqrt{b}\vert \leq \frac{\alpha}{2}. $$
Let $ K$ be a positive integer number such that $$\frac{{C} \sigma}{1+KC_{a,b}}< \delta.$$
Using Lemma \ref{Weyl}, \eqref{casek} and \eqref{max},
we get the following inequalities almost surely for all large $N$. 

\noindent 
If $i_N < n$,
\begin{eqnarray*}
\sqrt{\lambda _{i_N+1}(M_N)}&\leq &\sqrt{\lambda _{i_N+1}(M_N^{(K)})}
+(\sigma-\sigma _K)\sqrt{\lambda _{1}(X_NX_N^*)}\\
&<&\sqrt{s_K}+\frac{{C} \sigma}{1+KC_{a,b}}\\
&\leq & \sqrt{a}+\alpha
\end{eqnarray*}
Similarly, if $i_N > 0$, 
\begin{eqnarray*}
\sqrt{\lambda _{i_N}(M_N)}
&\geq & \sqrt{b}-\alpha.
\end{eqnarray*}
Since
almost surely, for all $N$ large enough, $[{a}; {b}]$ is a gap in the spectrum of $M_N$,
we can deduce that, almost surely, for all $N$ large enough 
\begin{equation}\label{inf}
\lambda _{i_N+1}(M_N) < {a} \mbox{~~~~if~~} i_N < n, 
\end{equation}
\begin{equation}\label{sup} 
\quad \text{and}\quad \lambda _{i_N}(M_N) > {b}\mbox{~~~~if~~} i_N > 0.
\end{equation}
Since $\lambda _{n+1}(M_N)=-\lambda _0(M_N)=-\infty $, 
\eqref{inf} (resp. \eqref{sup}) is obviously satisfied if $i_N=n$ (resp. $i_N=0$). 
This ends the proof of Proposition \ref{seppre}. \end{proof}

\noindent Now, if $c<1$ and $[a;b]$ satisfies assumption 7 of Theorem \ref{sep} and is such that $\omega_{\sigma,\nu,c}(a)\leq 0$ but $\omega_{\sigma,\nu,c}(b)>0$, 
the exact separation phenomenon  follows by applying Theorem \ref{pasde} to $[a;b]$ and  Proposition \ref{seppre} to $[a^{'};b]$ where $a<a^{'}<b$ is chosen such that 
$\omega_{\sigma,\nu,c}(a^{'})>0$. The proof of Theorem \ref{sep} is complete.\\

Note that for $c=1$, it turns out that Theorem \ref{pasde} and Theorem \ref{sep} can be extended to the case $a=0$ as follows.
\begin{proposition}\label{extc=1} Assume that $c=1$.
 Let $b>0 $ be such that 
\begin{itemize}
\item  there exists $\delta>0$ such that for all large $N$, $ [0; b+\delta[  \subset \mathbb{R}\setminus \rm{supp}( \mu_{\sigma,\mu _{A_N A_N^*},c_N})$
 \item   $A_{Nj}$ denoting the matrix resulting from removing the $j$-th column from $A_N$, there exists $0<\tau<\delta$ and a positive $d<1$ such that for all $N$ large, the number of j's with no eigenvalues 
of $N/(N-1) A_{Nj}A_{Nj}^*$ appearing in   $ \omega_{{\sigma,\mu _{A_N A_N^*},c_N}}([0, b +\tau[)$ is greater that $N-N^{d}$. 
\end{itemize}
Then  $[\omega_{\sigma,\nu,1}(0); \omega_{\sigma, \nu, 1}(b)]$ is to the left of the spectrum of $A_NA_N^*$ and 
 $$\mathbb P\left(\mbox{for all large N}, \rm{spect}(M_N) \subset  ]b;+\infty[ \right)=1.$$
\end{proposition}
\begin{proof} According to Lemma \ref{convg}, $[\omega_{\sigma,\nu,1}(0); \omega_{\sigma, \nu,1}(b+\frac{\delta}{2})] \subset    \mathbb{R}\setminus \rm{supp} 
(\mu_{A_NA_N^*})$ with $\omega_{\sigma,\nu,1}(0)=0$ and $ \omega_{\sigma,\nu,1}(b+\frac{\delta}{2}) >0$. 
Applying Proposition \ref{seppre} for $[b/2;b]$ yields that the spectrum  of $M_N$ is on the right hand side of $b$.\end{proof}


\noindent
\section{Application to spiked models}\label{spi}
\subsection{Matricial model and notations}
\noindent We will consider the deformed model:
\begin{equation}\nonumber 
M_N=( \sigma X_N+A_N)(\sigma X_N+A_N)^*\end{equation}
\begin{itemize}
\item $n \leq N$, $c_N=n/ N\rightarrow c \in]0;1].$
\item $X_N$ satisfies conditions 1., 2. and 3. of Theorem \ref{pasde}.
\item Let $\nu \neq \delta_0$ be a probability measure whose support is a finite union of disjoint (possibly degenerate) closed bounded intervals. Let $\theta _1 > \ldots > \theta _J\geq 0$ be
$J$ fixed real numbers 
independent of $N$ which are outside the support of $\nu $.
Let $k_1,\ldots,k_J$ be fixed  integer numbers independent of $N$ and  $r=\sum_{j=1}^J k_j$.
Let $\beta_j(N)\geq 0$, $r+1\leq j \leq n$, be such that
$\frac{1}{n} \sum_{j=r+1}^{n} \delta _{\beta _j(N)}$ weakly
converges to  $\nu $ and
 \begin{equation}\label{univconv}
 \max _{r+1\leq j\leq n} {\rm dist}(\beta _j(N),{\rm supp}(\nu ))\vers _{N \rightarrow \infty } 0.\end{equation}
Let $a_j(N), j=1,\ldots,J$, and  $b_j(N)\geq 0$, $r+1\leq j \leq n$, be complex numbers such that 
$$\forall j=1,\ldots, J, l=1,\ldots, k_j, \lim_{N \rightarrow +\infty} \vert a_{k_1+\cdots + k_{j-1} +l}(N)\vert^2 =\theta_j,$$
$$\vert b_j(N)\vert^2=\beta_j(N).$$
Let us introduce the $n \times N $ deterministic matrix $A_N$ by setting

  $$(A_N)_{ii}=a_i(N) \mbox{~~ for  $i=1,\ldots,r$},$$
 for any $r+1\leq i\leq n$,
$$(A_N)_{ii}=b_i(N).
$$ and else $(A_N)_{ij}=0.$

Thus, $A_NA_N^*$ is a diagonal matrix,
  each $\vert a_j(N) \vert^2$ 
is an eigenvalue of $A_N A_N^*$ with a fixed multiplicity $k_j$ (with $\sum_{j=1}^J k_j=r$) and the other eigenvalues of $A_N A_N^*$ are the 
$\beta_j(N)$. Moreover the empirical spectral measure of ${A_N A_N^*} $ weakly
converges to  $\nu $.
We will denote by $\Theta$ the set  $$\Theta:=\left\{ \theta _1 ; \cdots ; \theta _J \right\}$$
and by $K_{\nu,\Theta}$ the set  \begin{equation} \label{K} K_{\nu,\Theta}:= \mbox{supp}(\nu) \cup \Theta.\end{equation}
\end{itemize}
According to Proposition \ref{compsup}, there exists a nonnull integer number  $p$ and  $u_1< v_1<u_2<\ldots <u_p < v_p$  (depending on $\sigma,\nu,c$) such that 
$${{\cal E}_{\sigma,\nu,c}}=]-\infty;u_1[\cup_{l=1}^{p-1}]v_l;u_{l+1}[\cup 
]v_p;+\infty[,$$ $\mbox{supp}(\nu) \subset \cup_{l=1}^{p} [u_l;v_l]$
and for  each $l \in \{1,\ldots,p\}$,  $[u_l;v_l]\cap\mbox{supp}(\nu) \neq \emptyset$,
$$
 \mbox{supp}(\mu_{\sigma,\nu,c})
=\cup_{l=1}^p [\Phi_{\sigma,\nu,c}(u_l^-);\Phi_{\sigma,\nu,c}(v_l^+)],$$
with \\

$\Phi_{\sigma,\nu,c}(u_1^-) < \Phi_{\sigma,\nu,c}(v_1^+)< \Phi_{\sigma,\nu,c}(u_2^-) < \Phi_{\sigma,\nu,c}(v_2^+)$
$$\hspace*{6cm}
<\cdots< \Phi_{\sigma,\nu,c}(u_p^-) < \Phi_{\sigma,\nu,c}(v_p^+),$$
where $ \Phi_{\sigma,\nu,c}(u_l^-) =\lim_{u\uparrow u_l} \Phi_{\sigma,\nu,c}(u)$ and $ \Phi_{\sigma,\nu,c}(v_l^+) =\lim_{u\downarrow v_l} \Phi_{\sigma,\nu,c}(u)$.
\subsection{Subordination property of rectangular free convolution of ratio $c$}\label{convrect}
\addtocounter{footnote}{+1}
For any $c\in]0,1]$, the rectangular free convolution $\boxplus_{c}$ is defined in \cite{BG} 
in the following way. Let $M_{n,N}$, $N_{n,N}$ be $n$ by $N$ independent random matrices, one of  them having a distribution which is invariant by multiplication by any unitary matrix on any side, the symmetrized {empirical singular measures\footnotemark} of which tend, as $N$ tends to infinity in such a way that $n/N$ tends to $c$,
to nonrandom probability measures $\nu_1$, $\nu_2$.
\footnotetext{The empirical singular measure of a $n$ by $N$ $(n\leq N)$ matrix $M$ is the uniform law on the eigenvalues of $\sqrt{MM^*}$.}
   Then the  symmetrized  empirical singular law of $M_{n,N} +N_{n,N}$ tends to $\nu_1\boxplus_c \nu_2$.\\

For any probability measure  $\tau$ on $[0; +\infty[$, let us denote by $\sqrt{\tau}^s$ the symmetrization of the pushfoward of $\tau$ by the map $t \mapsto \sqrt{t}$ and by $\sigma^2 \tau$ the pushforward of $\tau$ by the map $t \mapsto \sigma^2 t.$
Note that the limiting spectral measure $\mu$ of $M_N$ defined in (\ref{modele}) does not depend on the distribution of the entries of $X_N$.
Thus, choosing gaussian entries for $X_N$, we can deduce that the limiting spectral measure $\mu_{\sigma,\nu,c}$ of $M_N$ satisfies 
$$\sqrt{\mu_{\sigma,\nu,c}}^s =\sqrt{\nu}^s \boxplus_c \sqrt{\sigma^2 \mu_c}^s,$$
\noindent where $\mu_c$ is the well-known
Marchenko-Pastur law defined by
\begin{equation} \label{MP} \mu_c(dx)=
f(x) \1_{[(1-\sqrt{c})^2;
(1+\sqrt{c})^2]}(x)dx\end{equation}
with $$ f(x)=\frac{\sqrt{\left(x-(1-\sqrt{c})^2\right)
\left((1+\sqrt{c})^2-x\right)}}{2\pi c x}.$$

Let us now explain the intuition of our result.
\noindent Using (\ref{TS3}), one can easily see that for any $x \in \mathbb{R}\setminus \{ \mbox{supp}({\mu_{\sigma,\nu,c}})\cup\{0\}\}$,
$$ c\omega_{\sigma,\nu,c}(x) g_\nu(\omega_{\sigma,\nu,c}(x))^2+(1-c)g_\nu(\omega_{\sigma,\nu,c}(x))=cx g_{\mu_{\sigma,\nu,c}}(x)^2 +(1-c) g_{\mu_{\sigma,\nu,c}}(x). $$
Note that this has to be related with the subordination result established in \cite{BBG}  involving the $H$ transform related to rectangular free convolution with ratio $c$. 
Indeed, for any probability measure  $\tau$ on $[0; +\infty[$, the $H$ transform of ${\sqrt{\tau}}^s$ is such that 
$H_{\sqrt{\tau}^s} (z) =\frac{c}{z} g_\tau(\frac{1}{z})^2 +(1-c) g_\tau(\frac{1}{z})$ and the authors established in  \cite{BBG} that for any probability measures $\mu$ and $\tau$ on  $[0; +\infty[$ such that $\tau$ is $\boxplus_c$ infinitely divisible, there exists a meromorphic function $F$ on $\mathbb{C}\setminus [0;+\infty[$ such that  
$H_{\sqrt{\mu}^s \boxplus_c \sqrt{\tau}^s}(z)=H_{\sqrt{\mu}^s}(F(z)).$

\noindent The intuition is that  for large $N$, \\

\noindent $cx g_{\mu_{M_N}}(x)^2 +(1-c) g_{\mu_{M_N}}(x) $
\begin{eqnarray*}& \approx &cx g_{\mu_{\sigma, \mu_{A_N A_N^*},c}}(x)^2 +(1-c) g_{\mu_{\sigma, \mu_{A_N A_N^*},c}}(x) \\&\approx & 
 c\omega_{\sigma,\nu,c}(x) g_{\mu _{A_N A_N^*}}(\omega_{\sigma,\nu,c}(x))^2+(1-c)g_{\mu _{A_N A_N^*}}(\omega_{\sigma,\nu,c}(x))
\end{eqnarray*}
and therefore that there will be  eigenvalues of $M_N$  that separate from the bulk whenever some of the  equations $$
\omega_{\sigma,\nu,c}(x)=\theta_j,$$ $j=1,\ldots,J$, admits a  solution outside the support of $\mu_{\sigma,\nu,c}$. According to Theorem \ref{caractfinale}, $\omega_{\sigma,\nu,c}(x)=\theta_j$ admits a solution outside the support of $\mu_{\sigma,\nu,c}$ if and only if 
$ \Phi_{\nu,\sigma,c}^{'}(\theta_j) >0, g_{\nu}(\theta_j) >-\frac{1}{\sigma^2c}.$ Moreover there is only one such a solution which is $\Phi_{\sigma,\nu,c}(\theta_j)$.\\
Thus, as in the case of additive or multiplicative  deformed square  models (see \cite{CDFF, C,  BBCF}), free subordination property definitely sheds light on the problem of outliers for information-plus-noise type models since, denoting by $F$ the rectangular $\boxplus_c$- free subordination function relative to $\nu$ and $\mu_{c}$, 
for large $N$, the $\theta_i$'s such that  
the equation $$\frac{1}{F(\frac{1}{\rho})}=\theta_i$$ has solution $\rho$ outside  the support of $\mu_{\sigma,\nu,c}$ will actually  generate $k_i$  eigenvalues of $M_N$  in  a neighborhood of this $\rho$.
\noindent This intuition leads us to introduce the following set
$$\Theta_{\sigma,\nu,c} = \left\{\theta \in \Theta,   \Phi_{\sigma,\nu,c}^{'}(\theta) >0, g_\nu(\theta) >-\frac{1}{\sigma^2c}\right\},$$
\noindent and to introduce  for any $\theta$ in $\Theta_{\sigma,\nu,c}$, $$\rho_{\theta} =\Phi_{\sigma,\nu,c}(\theta).$$
\noindent According to B)  of Theorem \ref{caractfinale},  $\forall \theta \in \Theta_{\sigma,\nu,c}$, $\rho_{\theta}\notin \mbox{supp}(\mu_{\sigma,\nu,c}).$
\subsection{Inclusion of the spectrum}
Let us  define \begin{equation}\label{defSsigma}{\cal S}=\mbox{~supp~} (\mu_{\sigma,\nu,c}) \cup \left\{ \rho_{\theta}, \theta \in \Theta_{\sigma,\nu,c} \right\}.\end{equation}
In order to establish an inclusion of the spectrum of $M_N$,  in the following Theorem \ref{inclusion}, 
 we first prove the next lemma.
\begin{lemme}\label{inclusionN}
Let $[u,v] \subset   \mathbb{R}\setminus {\cal S}$. 
For any  $\delta>0$ small enough,  for all large $N$, $[u-\delta; v+\delta] \subset  \mathbb{R}\setminus   \mbox{supp}(\mu_{\sigma, \mu_{A_N A_N^*},c_N}) $.
\end{lemme}

\begin{proof}
Let  $\delta_0  >0$ be such that $[u-\delta_0;v+\delta_0] \subset   \mathbb{R}\setminus {\cal S}$.  According to B)  of  Theorem \ref{caractfinale}, we have $]\omega_{\sigma, \nu,c}(u-\delta_0);\omega_{\sigma,\nu,c}(v+\delta_0)[\subset  \mathbb{R}\setminus  \mbox{supp}(\nu)$ and for any $x$ in $]\omega_{\sigma,\nu,c}(u-\delta_0);\omega_{\sigma,\nu,c}(v+\delta_0)[$, $\Phi_{\sigma,\nu,c}^{'}(x) >0$, $g_\nu(x) >-\frac{1}{\sigma^2 c}$. In particular $\left\{\Theta \setminus \Theta_{\sigma,\nu,c}\right\} \cap ]\omega_{\sigma,\nu,c}(u-\delta_0);\omega_{\sigma,\nu,c}(v+\delta_0)[ =\emptyset$.
Moreover
 it is clear by  B)  of  Theorem \ref{caractfinale} that, since for any $\theta$ in $\Theta_{\sigma,\nu,c}$, $\rho_{\theta} \notin ]u-\delta_0;v+\delta_0[$, we have
$ \Theta_{\sigma,\nu,c} \cap ]\omega_{\sigma,\nu,c}(u-\delta_0);\omega_{\sigma,\nu,c}(v+\delta_0)[ =\emptyset$. Therefore, $$ ]\omega_{\sigma,\nu,c}(u-\delta_0);\omega_{\sigma,\nu,c}(v+\delta_0)[ \subset   \mathbb{R}\setminus K_{\nu, \Theta},$$ \noindent where $ K_{\nu, \Theta}$ is defined in (\ref{K}).
\noindent Note that using  B)  of  Theorem \ref{caractfinale}, we have 
$ [\omega_{\sigma,\nu,c}(u-\frac{\delta_0}{2});\omega_{\sigma,\nu,c}(v+\frac{\delta_0}{2})] \subset  ]\omega_{\sigma,\nu,c}(u-\delta_0);\omega_{\sigma,\nu,c}(v+\delta_0)[.$ Then there exists $\alpha>0$ such that $$d\left(  [\omega_{\sigma,\nu,c}(u-\frac{\delta_0}{2});\omega_{\sigma,\nu,c}(v+\frac{\delta_0}{2})],  K_{\nu, \Theta}\right) > \alpha.$$
According to the assumptions on the eigenvalues of $A_N A_N^*$, for all large $N$, we have 
$$\max \{\max_{r+1\leq i \leq n} d(\beta_i(N), \rm{supp}(\nu)); \max_{j=1,\ldots,J} \max_{l=1,\ldots,k_j}d(\vert a_{k_1+\cdots+ k_{j-1}+l}(N) \vert^2 , \theta_j)\}< \alpha/2$$ so that 
the spectrum of  $A_N A_N^*$ is included in $\left\{x, d\left( x,  K_{\nu, \Theta}\right)< \frac{\alpha}{2}\right\}.$
Moreover there exists $\epsilon>0$ such that \begin{equation}\label{strict} \mbox{ for any $x$ in $ [\omega_{\sigma,\nu,c}(u-\frac{\delta_0}{2});\omega_{\sigma,\nu,c}(v+\frac{\delta_0}{2})]$,}~~  g_\nu(x) >-\frac{1}{\sigma^2 c}+ \epsilon \mbox{    and~~} \Phi_{\sigma,\nu,c}^{'}(x) > \epsilon.\end{equation}

\noindent Since,  for all large $N$, $g_{\mu _{A_NA_N^*}}$ 
 is analytic and  uniformly bounded  on $\{ z\in \C , {\rm dist}(z,K_{\nu,\Theta}) > {\alpha }\}$, Vitali's and Weierstrass' theorems yield that 
 $g_{\mu_{A_N A_N^*}}$ and $g_{\mu_{A_N A_N^*}}^{'}$  converge to $g_\nu$ and   $g_\nu '$, respectively, uniformly on $ [\omega_{\sigma,\nu,c}(u-\frac{\delta_0}{2});\omega_{\sigma,\nu,c}(v+\frac{\delta_0}{2})]$.
The uniform convergence of  $ \Phi^{'}_{\sigma,\mu_{A_N A_N^*},c_N}$ 
towards $\Phi_{\sigma,\nu,c}^{'}$  on $ [\omega_{\sigma,\nu,c}(u-\frac{\delta_0}{2});\omega_{\sigma,\nu,c}(v+\frac{\delta_0}{2})]$ readily follows. 
 Hence, using (\ref{strict}),
 we can claim that for all large $N$, for all $x$ in  $ [\omega_{\sigma,\nu,c}(u-\frac{\delta_0}{2});\omega_{\sigma,\nu,c}(v+\frac{\delta_0}{2})]$,
$g_{\mu_{A_N A_N^*}}(x) >-\frac{1}{\sigma^2 c}+ \frac{\epsilon}{2}>-\frac{1}{\sigma^2 c_N}+ \frac{\epsilon}{4} \mbox{    and~~} \Phi^{'}_{\sigma,\mu_{A_N A_N^*},c_N}(x) > \frac{\epsilon}{2}$.
Therefore for all large $N$, \\

$ [\omega_{\sigma,\nu,c}(u-\frac{\delta_0}{2});\omega_{\sigma,\nu,c}(v+\frac{\delta_0}{2})]$
$$\subset \left\{ u \in   \mathbb{R}\setminus \mbox{supp}(\mu _{A_N A_N^*}),  \Phi_{\sigma,\mu _{A_N A_N^*}, c_N}^{'}(u) >0, g_{\mu _{A_N A_N^*}}(u) >-\frac{1}{\sigma^2 c_N}\right\}.$$
According to   B)  of  Theorem \ref{caractfinale}, we can deduce that
\begin{eqnarray*}&&\Phi_{\sigma, \mu _{A_N A_N^*},c_N}\left( [\omega_{\sigma,\nu,c}(u-\frac{\delta_0}{2});\omega_{\sigma,\nu,c}(v+\frac{\delta_0}{2})]\right)\\ &&\hspace*{1.5cm}=\left[ [\Phi_{\sigma,\mu _{A_N A_N^*},c_N}(\omega_{\sigma,\nu,c}(u-\frac{\delta_0}{2}));\Phi_{\sigma,\mu _{A_N A_N^*},c_N}(\omega_{\sigma,\nu,c}(v+\frac{\delta_0}{2}))\right]
\\&&\hspace*{2cm}\subset   \mathbb{R}\setminus  \mbox{supp}(\mu_{\sigma,\mu_{A_NA_N^*},c_N}).\end{eqnarray*}
Now since $\Phi_{\sigma,\mu _{A_N A_N^*},c_N}(\omega_{\sigma,\nu,c}(u-\frac{\delta_0}{2}))$ and $\Phi_{\sigma,\mu _{A_N A_N^*},c_N}(\omega_{\sigma,\nu,c}(v+\frac{\delta_0}{2}))$ converge respectively towards $\Phi_{\sigma,\nu,c}(\omega_{\sigma,\nu,c}(u-\frac{\delta_0}{2}))$ and $\Phi_{\sigma,\nu,c}(\omega_{\sigma,\nu,c}(v+\frac{\delta_0}{2}))$ and using B) of Theorem \ref{caractfinale}, we have 
for all large $N$, $$\Phi_{\sigma,\mu _{A_N A_N^*},c_N}(\omega_{\sigma,\nu,c}(u-\frac{\delta_0}{2}))\leq u-\frac{\delta_0}{4},
$$ 
$$\Phi_{\sigma,\mu _{A_N A_N^*},c_N}(\omega_{\sigma,\nu,c}(v+\frac{\delta_0}{2})) \geq v+\frac{\delta_0}{4}$$
\noindent and then $$[u-\frac{\delta_0}{4};v+\frac{\delta_0}{4} ] \subset      \mathbb{R}\setminus  \mbox{supp}(\mu_{\sigma,\mu_{A_NA_N^*},c_N})$$
\noindent and the proof of Lemma \ref{inclusionN} is complete.\end{proof}
We have the following inclusion of the spectrum of $M_N$.
\begin{theoreme}\label{inclusion}
A) For any  $\epsilon>0$, 
$$\mathbb P[\mbox{for all large N}, \rm{spect}(M_N) \subset \{x \in \R , \mbox{dist}(x,{\cal S}\cup \{0\})\leq \epsilon \}]=1.$$
B) If $c<1$ and if moreover  $\theta_J>0$ and   $0 \in {\cal E}_{\sigma,\nu,c}$, then,  for any  $\epsilon>0$, $$\mathbb P[\mbox{for all large N}, \rm{spect}(M_N) \subset \{x \in \R , \mbox{dist}(x,{\cal S})\leq \epsilon \}]=1.$$
C) If $c=1$, $$\mathbb P[\mbox{for all large N}, \rm{spect}(M_N) \subset \{x \in \R , \mbox{dist}(x,{\cal S})\leq \epsilon \}]=1.$$
\end{theoreme}
\begin{proof} 
Using (\ref{max}), assumption 4. of Theorem \ref{pasde} and Lemma \ref{Weyl},  one can easily see that
there exists $R>0$ such that 
${\cal S} \subset [0;R[$ and 
almost surely for all large $N$, $\Vert M_N \Vert < R$. Let us fix $\epsilon>0$ such that 
 $$ 2 \epsilon  < \min\left \{\min_{l=1,\ldots, p-1} [\Phi_{\sigma,\nu,c}(u_{l+1}^-)-\Phi_{\sigma,\nu,c}(v_l^+)], d(\rho_{\theta_i}, \mbox{supp~} (\mu_{\sigma,\nu,c})),\right.$$
 \begin{equation}\label{eps} \hspace*{5cm} d(\rho_{\theta_i},\rho_{\theta_j}),\theta_i \neq \theta_j \mbox{~in~} 
\Theta_{\sigma,\nu,c}, R-\max {\cal S}\Big\}\end{equation}
If  $\min{\cal S}>0$, choose  $\epsilon$ such that (\ref{eps}) holds and moreover $2 \epsilon < \min{\cal S}$.

\begin{itemize}
\item There exists $0\leq y_1<x_1<\cdots< y_q<x_q \leq \max {\cal S}$, such that 
$$[0; R]\setminus \{x \in \R , \mbox{dist}(x,{\cal S}\cup \{0\})< \epsilon \}=   \cup_{i=1}^{q} [y_i+\epsilon; x_i -\epsilon] \cup [\max \{u \in {\cal S}\} +\epsilon;R].$$
Applying Theorem \ref{pasde} for each of these intervals (since  Lemma \ref{inclusionN} implies 7. (i) of Theorem \ref{pasde} and 7. (ii) of Theorem \ref{pasde} is obviously satisfied since $A_N$ is of the form (\ref{diagonale})), we get A).
\item
Assume that $c=1$. If $0 \in {\cal S}$, A) implies C). If $\min {\cal S}>0$, we can moreover apply Proposition \ref{extc=1} to $[0;\min {\cal S}-\epsilon]$ and,
using also A), deduce C).
\item Now, assume that $0<c<1$. If $0 \in {\cal E}_{\sigma,\nu,c}$ (or equivalently $u_1>0$ where $u_1$ is defined in Theorem \ref{caractfinale} D))   and if $\theta_J >0$,
let $u=\min \{ u_1, \theta_J \}>0.$
Applying   Lemma \ref{inclusionN}  and Theorem \ref{sep} to $[ \Phi_{\sigma,\nu,c}(\frac{u}{4}); \Phi_{\sigma,\nu,c}(\frac{3u}{4})]$, we obtain that almost surely for all large $N$, there is no eigenvalue of $M_N$ on the left hand side  of $\Phi_{\sigma,\nu,c}(\frac{u}{4})$ since $\mbox{supp}(\nu)\subset\cup_{l=1}^p[u_l;v_l]$ and using the assumptions 
on the spectrum of $A_NA_N^*$.
Using also A), we deduce B).
\end{itemize} \end{proof}

\subsection{Convergence of eigenvalues }
\noindent 
In the non-spiked case  i.e. $r=0$, 
the result of Theorem  \ref{inclusion} reads as: \begin{itemize}
\item If $c<1$, \\
$\forall \epsilon > 0$, 
\begin{eqnarray}{\label{InclNonPerturb}}
\mathbb P[{\rm Spect}(M_N)\subset {\rm supp}(\mu _{\sigma,\nu,c} )\cup\{0\}+(-\epsilon ,\epsilon ), 
\, \text{for all $N$ large}]=1,
\end{eqnarray}
If $u_1>0$, then 
\begin{eqnarray}{\label{InclNonPerturb2}}
\mathbb P[{\rm Spect}(M_N)\subset {\rm supp}(\mu _{\sigma,\nu,c} )+(-\epsilon ,\epsilon ), 
\, \text{for all $N$ large}]=1.
\end{eqnarray}
\item If $c=1$, then 
\begin{eqnarray}{\label{InclNonPerturb2}}
\mathbb P[{\rm Spect}(M_N)\subset {\rm supp}(\mu _{\sigma,\nu,c} )+(-\epsilon ,\epsilon ), 
\, \text{for all $N$ large}]=1.
\end{eqnarray}
\end{itemize}
This readily leads to the following asymptotic result for the extremal eigenvalues.

\begin{proposition}{\label{ThmASCVNonSpike}}
Assume that the deformed model $M_N$ is without spike i.e. $r=0$. Let $k\geq 0$ be a fixed integer.
\begin{enumerate}
\item The first largest eigenvalues $\lambda_{1+k}(M_N)$ 
 converge almost surely to the right
endpoint of the support of $\mu _{\sigma, \nu,c} $.
\item If $\Phi_{\sigma , \nu,c }(u_1^-)=0$ that is when   $c=1$
and either $0\in \rm{supp}( \nu)$ or $0\notin \rm{supp}( \nu)$ and  $g_\nu(0) \leq -\frac{1}{\sigma^2}$, 
then the last smallest eigenvalues $\lambda_{n-k} (M_N)$  converge almost surely to zero.
\item If $u_1>0$ (which implies $\Phi_{\sigma,\nu,c}( u_1^-) >0$) then the last  smallest eigenvalues $\lambda_{n-k} (M_N)$  converge almost surely to $\Phi_{\sigma , \nu,c }(u_1^-)$.
\end{enumerate}
\end{proposition}

\begin{proof} 
Recalling that ${\rm supp}(\mu _{\sigma, \nu,c} )=\cup _{l=1}^p[\Phi_{\sigma , \nu,c }(u_l^-),\Phi_{\sigma , \nu,c }(v_l^+)]$, 
from \eqref{InclNonPerturb}, one has that, for all $\epsilon > 0$, 
$$\mathbb P[\limsup _N\lambda _1(M_N)\leq \Phi_{\sigma , \nu }(v_p^+)+\epsilon ]=1.$$
But as $\Phi_{\sigma , \nu,c }(v_p^+)$ is a boundary point of ${\rm supp}(\mu _{\sigma, \nu,c} )$, 
the number of eigenvalues of $M_N$ falling into $[\Phi_{\sigma , \nu,c }(v_p^+)-\epsilon ,\Phi_{\sigma , \nu,c }(v_p^+)+\epsilon ]$ 
tends almost surely to infinity as $N\to \infty $. Thus, almost surely, 
$$\liminf _N\lambda _{1+k}(M_N)\geq \Phi_{\sigma , \nu,c }(v_p^+)-\epsilon .$$
1)  follows by letting $\epsilon \to 0$. \\
 The proofs of 2) and 3) are  similar to the proof of 1) using the fact that in these cases we have
\begin{eqnarray*}
\mathbb P[{\rm Spect}(M_N)\subset {\rm supp}(\mu _{\sigma,\nu,c} )+(-\epsilon ,\epsilon ), 
\, \text{for all $N$ large}]=1.
\end{eqnarray*}
\end{proof}

In the spiked case, in order 
to obtain, in the following  Theorem \ref{ThmASCV2}, a  description of the convergence of the
eigenvalues of $M_N$, depending on the location of the spikes of the perturbation, we first note that one can  readily deduce the following corollary from Theorem \ref{sep} and Lemma \ref{inclusionN}.

\begin{corollaire}\label{thetat}
Let us fix $\epsilon>0$ such that 
 $$ 2 \epsilon  < \min\left \{\min_{l=1,\ldots, p-1} [\Phi_{\sigma,\nu,c}(u_{l+1}^-)-\Phi_{\sigma,\nu,c}(v_l^+)], d(\rho_{\theta_i}, \mbox{supp~} (\mu_{\sigma,\nu,c})),\right.$$
 \begin{equation}\nonumber \hspace*{7cm} d(\rho_{\theta_i},\rho_{\theta_j}),\theta_i \neq \theta_j \mbox{~in~} 
\Theta_{\sigma,\nu,c}\Big\}\end{equation} Let $u$ be  in $\Theta_{\sigma,\nu,c}\cup \{ v_l, l=1, \ldots , m\} $ 
(resp.  if $c =1$ in  $\Theta_{\sigma,\nu,c}\cup \{ u_l, l=1, \ldots , m\} $ and if $c<1$ in $\left(\Theta_{\sigma,\nu,c}\cup \{ u_l, l=1, \ldots , m\} \right)\cap ]0;+\infty[$).
Let us choose $\delta > 0$ small enough so that for large $N$, 
$[u+\delta ; u+2\delta ]$ (resp. $[u-2\delta ; u-\delta ]$) 
is included  in ${\cal E}_{\sigma,\nu,c}$ if $c=1$ and  in ${\cal E}_{\sigma,\nu,c}\cap ]0;+\infty[$ if $c<1$  
and for any $0\leq \delta '\leq 2\delta $, $\Phi_{\sigma , \nu,c }(u+\delta ')-\Phi_{\sigma , \nu,c }(u^+) < \epsilon $ 
(resp. $\Phi_{\sigma , \nu,c }(u^-)-\Phi_{\sigma , \nu,c }(u-\delta ') < \epsilon $). 
For $N$ large enough, let $i_N=i_N(u)$ be such that 
$$\lambda _{i_N+1}(A_NA_N^*) < u+\delta \, \text{ and } \, \lambda _{i_N}(A_NA_N^*) > u+2\delta $$
(resp. $\lambda_{i_N+1}(A_NA_N^*) < u-2\delta \, \text{ and } \, \lambda _{i_N}(A_NA_N^*) > u-\delta $). 
Then $$\mathbb P\big{[}\lambda_{i_N+1}(M_N) < \Phi_{\sigma , \nu,c }(u+\delta) \, \text{ and } \, 
\lambda _{i_N}(M_N) > \Phi_{\sigma , \nu,c }(u+2 \delta), \text{ for all large $N$}\big{]}=1.$$
(resp. $\mathbb P\big{[}\lambda _{i_N+1}(M_N) < \Phi_{\sigma , \nu,c }(u-2\delta) 
 \, \text{ and } \, 
\lambda _{i_N}(M_N) > \Phi_{\sigma , \nu,c }(u-\delta)  \text{ for all large $N$}\big{]}=1.$)
\end{corollaire}
\begin{theoreme}{\label{ThmASCV2}}
For any $j=1,\ldots,J$,
we denote by $n_{j-1}+1, \ldots , n_{j-1}+k_j$ the descending ranks of $\{\vert  a_{k_1+\cdots+k_{j-1}+l}(N)\vert^2, l=1,\ldots,k_j\}$ among the eigenvalues of $A_N A_N^*$.
\begin{itemize}
\item[{1)}] If $\theta_j \in {\cal E}_{\sigma,\nu,c}$ 
(i.e. $\in \Theta _{\sigma,\nu,c}  $), and moreover, if $c<1$ $\theta_j \neq 0$, then the $k_j$ eigenvalues $(\lambda_{n_{j-1}+i}(M_N), \, 1 \leq i \leq k_j)$ 
converge almost surely outside the support of $\mu _{\sigma,\nu,c} $ 
towards $\rho _{\theta _j}=\Phi_{\sigma , \nu,c }(\theta _j)$.
\item[\text{2)}] If $\theta_j \in \mathbb{R} \setminus  {\cal E}_{\sigma,\nu,c}$, 
then we let $[u_{l_j}, v_{l_j}]$ (with $1\leq l_j\leq p$) be the connected component 
of $ \mathbb{R} \setminus  {\cal E}_{\sigma,\nu,c}$ which contains $\theta _j$.
\begin{itemize}
\item[{a)}] If $\theta_j$ is on the right of any connected component of ${\rm supp}(\nu )$ 
which is included in $[u_{l_j},v_{l_j}]$ then the $k_j$ eigenvalues $(\lambda_{n_{j-1}+i}(M_N)$, $1\leq i\leq k_j)$ 
converge almost surely to $\Phi_{\sigma , \nu,c }(v_{l_j}^-)$ 
which is a boundary point of the support of $\mu _{\sigma, \nu,c} $. \\
If $u_{l_j}>0$ (which is always true if $ l_j\neq 1$)  and if $\theta_j$ is on the left of any connected component of ${\rm supp}(\nu )$ 
which is included in $[u_{l_j},v_{l_j}]$ then the $k_j$ eigenvalues $(\lambda_{n_{j-1}+i}(M_N)$, $1\leq i\leq k_j)$ 
converge almost surely to $\Phi_{\sigma , \nu,c }(u_{l_j}^{-})$ which is a boundary point of the support of $\mu _{\sigma, \nu,c}$.
\item[{b)}]If $ l_j =1$ and  $\Phi_{\sigma , \nu,c }(u_{1}^{-})=0$  and if $\theta_j$ is on the left of any connected component of ${\rm supp}(\nu )$ 
which is included in $[u_{1},v_{1}]$ then  the $k_j$ eigenvalues $(\lambda_{n_{j-1}+i}(M_N)$, $1\leq i\leq k_j)$ 
converge almost surely to 0.
\item[{c)}] If $\theta_j$ is between two connected components of ${\rm supp}(\nu )$ 
which are included in $[u_{l_j},v_{l_j}]$ then 
the $k_j$ eigenvalues $(\lambda_{n_{j-1}+i}(M_N)$, $1\leq i\leq k_j)$ 
converge almost surely to the $\alpha _j$-th quantile of $\mu _{\sigma, \nu,c} $ 
(that is to $q_{\alpha _j}$ defined by $\alpha _j=\mu _{\sigma, \nu,c} (]-\infty , q_{\alpha _j}])$) 
where $\alpha _j$ is such that $\alpha _j=1-\lim _N\frac{n_{j-1}}{N}=\nu (]-\infty ,\theta _j])$.
\end{itemize}
\end{itemize}
\end{theoreme}
\begin{proof}
 The proof follows the lines of the proof of Theorem 8.1 \cite{CDFF}. We include it for the reader's convenience.\\
%
\noindent 1) Choosing $u=\theta _j$  in Corollary \ref{thetat} gives, for any $\epsilon > 0$, 
\begin{eqnarray}\label{eqfin}
\rho _{\theta _j}-\epsilon \leq \lambda _{n_{j-1}+k_j}(M_N)\leq \cdots 
\leq \lambda _{n_{j-1}+1}(M_N)\leq \rho _{\theta _j}+\epsilon , \text{ for large $N$}
\end{eqnarray}
holds almost surely. Hence
$$\forall 1\leq i\leq k_j, \quad \lambda _{n_{j-1}+i}(M_N) \overset{a.s.}{\longrightarrow }\rho _{\theta _j}.$$

\noindent 
2) a) We only focus on the case where $\theta _j$ is on the right of any connected component 
of ${\rm supp}(\nu )$ which is included in $[u_{l_j}, v_{l_j}]$ 
since the other case may be considered with similar arguments. 
Let us consider the set $\{ \theta _{j_0} > \ldots > \theta_{j_p}\}$ of all the $\theta _i$'s 
being in $[u_{l_j}, v_{l_j}]$ and on the right of any connected component of ${\rm supp}(\nu )$ 
which is included in $[u_{l_j}, v_{l_j}]$. 
Note that we have for all large $N$, for any $0\leq h\leq p$, 
$$n_{j_h-1}+k_{j_h} =n_{j_h}$$ 
\noindent 
and $\theta _{j_0}$ is the largest eigenvalue of $A_NA_N^*$ which is lower than $v_{l_j}$.
\noindent 
 Applying Corollary \ref{thetat} with $u=v_{l_j}$, we get that, almost surely, 
$$\lambda _{n_{j_0-1}+1}(M_N) < \Phi_{\sigma , \nu,c }(v_{l_j}+\delta) \text{ and } 
\lambda _{n_{j_0-1}}(M_N) > \Phi_{\sigma , \nu,c }(v_{l_j}+2 \delta)\text{ for all large $N$.}$$ 
Now, almost surely, the number of eigenvalues of $M_N$ being in 
$]\Phi_{\sigma , \nu,c }(v_{l_j}^+)-\epsilon , \Phi_{\sigma , \nu,c }(v_{l_j}+\delta)]$ 
should tend to infinity when $N$ goes to infinity. 
Since almost surely for all large $N$, $\lambda _{n_{j_0-1}}(M_N) > \Phi_{\sigma , \nu,c }(v_{l_j}+2 \delta)$ and 
$\lambda _{n_{j_0-1}+1}(M_N) < \Phi_{\sigma , \nu,c }(v_{l_j}+\delta) $, we should have 
$$\hspace*{-2cm}\Phi_{\sigma , \nu,c }(v_{l_j}^+)-\epsilon \leq \lambda _{n_{j_p-1}+k_{j_p}}(M_N)\leq \ldots 
\leq \lambda _{n_{j_0-1}+1}(M_N)$$
$$\hspace*{7cm} < \Phi_{\sigma , \nu,c }(v_{l_j}+\delta)< \Phi_{\sigma , \nu,c }(v_{l_j}^+)+\epsilon.$$
Hence, we deduce that: $\forall 0\leq h\leq p$ and $\forall 1\leq i\leq k_{j_h}$, 
$ \lambda _{n_{j_h-1}+i}(M_N) \overset{a.s.}{\longrightarrow }\Phi_{\sigma , \nu,c }(v_{l_j}^+).$ 
The result then follows since $j \in \{ j_0,\ldots,j_p\} $. \\
\noindent b)  In this case, $\theta_i$ is one of the finite number of lowest eigenvalues of $A_N$. Then b) readily follows from the fact that  the number of eigenvalues of $M_N$ being in 
$[0 ,\epsilon]$ 
should tend to infinity when $N$ goes to infinity. \\
\noindent 
c) Let $\alpha _j=1-\lim _N\frac{n_{j-1}}{N}=\nu (]-\infty , \theta _j])$. 
Denote by $Q$ (resp. $Q_N$) the distribution function of $\mu _{\sigma ,\nu,c} $ 
(resp. of the spectral measure of $M_N$). 
$Q$ is continuous on $\R$ and strictly increasing on 
each interval $[\Phi _{\sigma , \nu,c }(u_l^-), \Phi _{\sigma , \nu,c }(v_l^+)], 1\leq l\leq m$.\\
From (\ref{palier}) and the hypothesis on $\theta_j$, 
$\alpha _j \in ]Q(\Phi _{\sigma , \nu,c }(u_{l_j}^-)), Q(\Phi _{\sigma , \nu,c }(v_{l_j}^+))[$ 
and there exists a unique $q_j\in ]\Phi _{\sigma , \nu,c }(u_{l_j}^-), \Phi _{\sigma , \nu,c }(v_{l_j}^+)[$ 
such that $Q(q_j)=\alpha _j$. Moreover, $Q$ is strictly increasing in a neighborhood of $q_i$. \\
Let $\epsilon > 0$. From the almost sure convergence of $\mu _{M_N}$ to $\mu _{\sigma ,\nu,c} $ , 
we deduce $$Q_N(q_j+\epsilon )\vers _{N\rightarrow \infty }Q(q_j+\epsilon ) > \alpha _j, \quad \text{a.s.}.$$
From the definition of $\alpha_j$, it follows that for large $N$, 
$N, N-1, \ldots , n_{j-1}+k_j, \ldots , n_{j-1}+1$ belong to the set $\{ k, \lambda _k(M_n)\leq q_j+\epsilon \}$ 
and thus, $$\limsup _{N\vers \infty }\lambda _{n_{j-1}+1}(M_N)\leq q_j+\epsilon .$$
In the same way, since $Q_N(q_j-\epsilon )\vers _{N\rightarrow \infty }Q(q_j-\epsilon ) < \alpha _j$,
$$\liminf _{N\vers \infty }\lambda _{n_{j-1}+k_j}(M_N)\geq q_j-\epsilon .$$
Thus, the $k_j$ eigenvalues $(\lambda_{n_{j-1}+i}(M_N)$, $1\leq i\leq k_j)$ 
converge almost surely to $q_j$. \end{proof}

\noindent {\bf Acknowledgments} We would like to thank the anonymous referee  for his  pertinent comments which led to an overall improvement of the paper.

\bibliographystyle{abbrv}
\bibliography{5432}

\begin{thebibliography}{10}

\bibitem{BaiSil99}
Z.~Bai and J.~W. Silverstein.
\newblock Exact separation of eigenvalues of large-dimensional sample
  covariance matrices.
\newblock {\em Ann. Probab.}, 27(3):1536--1555, 1999.

\bibitem{BaiSilverbook}
Z.~Bai and J.~W. Silverstein.
\newblock {\em Spectral {A}nalysis of {L}arge {D}imensional {R}andom
  {M}atrices}.
\newblock Mathematics {M}onograph {S}eries 2. Science Press, Beijing, 2006.

\bibitem{BaiSilver}
Z.~Bai and J.~W. Silverstein.
\newblock No eigenvalues outside the support of the limiting spectral
  distribution of information-plus-noise type matrices.
\newblock {\em Random Matrices Theory Appl.}, 1(1):1150004, 44, 2012.

\bibitem{BBG}
S.~Belinschi, F.~Benaych-Georges, and A.~Guionnet.
\newblock Regularization by free additive convolution, square and rectangular
  cases.
\newblock {\em Complex Anal. Oper. Theory}, 3(3):611--660, 2009.

\bibitem{BBCF}
S.~Belinschi, H.~Bercovici, M.~Capitaine, and M.~F{\'e}vrier.
\newblock Outliers in the spectrum of large deformed unitarily invariant
  models.
\newblock preprint \texttt{arXiv:1207.5443}, 2012.

\bibitem{BG}
F.~Benaych-Georges.
\newblock Rectangular random matrices, related convolution.
\newblock {\em Probab. Theory Related Fields}, 144(3-4):471--515, 2009.

\bibitem{BR}
F.~Benaych-Georges and R.~R. Nadakuditi.
\newblock The singular values and vectors of low rank perturbations of large
  rectangular random matrices.
\newblock {\em J. Multivariate Anal.}, 111:120--135, 2012.

\bibitem{C}
M.~Capitaine.
\newblock Additive/multiplicative free subordination property and limiting
  eigenvectors of spiked additive deformations of {W}igner matrices and spiked
  sample covariance matrices.
\newblock {\em J. Theoret. Probab.}, 26(3):595--648, 2013.

\bibitem{CDM}
M.~Capitaine and C.~Donati-Martin.
\newblock Free {W}ishart processes.
\newblock {\em J. Theoret. Probab.}, 18(2):413--438, 2005.

\bibitem{CDF09}
M.~Capitaine, C.~Donati-Martin, and D.~F{\'e}ral.
\newblock The largest eigenvalues of finite rank deformation of large {W}igner
  matrices: convergence and nonuniversality of the fluctuations.
\newblock {\em Ann. Probab.}, 37(1):1--47, 2009.

\bibitem{CDFF}
M.~Capitaine, C.~Donati-Martin, D.~F{\'e}ral, and M.~F{\'e}vrier.
\newblock Free convolution with a semicircular distribution and eigenvalues of
  spiked deformations of {W}igner matrices.
\newblock {\em Electron. J. Probab.}, 16:no. 64, 1750--1792, 2011.

\bibitem{CSBD}
R.~Couillet, J.~W. Silverstein, Z.~Bai, and M.~Debbah.
\newblock Eigen-inference for energy estimation of multiple sources.
\newblock {\em IEEE Trans. Inform. Theory}, 57(4):2420--2439, 2011.

\bibitem{DozierSilver2}
R.~B. Dozier and J.~W. Silverstein.
\newblock Analysis of the limiting spectral distribution of large dimensional
  information-plus-noise type matrices.
\newblock {\em J. Multiv. Anal.}, 98(6):1099--1122, 2007.

\bibitem{DozierSilver}
R.~B. Dozier and J.~W. Silverstein.
\newblock On the empirical distribution of eigenvalues of large dimensional
  information-plus-noise-type matrices.
\newblock {\em J. Multiv. Anal.}, 98(4):678--694, 2007.

\bibitem{HLN}
W.~Hachem, P.~Loubaton, and J.~Najim.
\newblock Deterministic equivalents for certain functionals of large random
  matrices.
\newblock {\em Ann. Appl. Probab.}, 17(3):875--930, 2007.

\bibitem{LV}
P.~Loubaton and P.~Vallet.
\newblock Almost sure localization of the eigenvalues in a {G}aussian
  information plus noise model---application to the spiked models.
\newblock {\em Electron. J. Probab.}, 16:no. 70, 1934--1959, 2011.

\bibitem{VLM}
P.~Vallet, P.~Loubaton, and X.~Mestre.
\newblock Improved subspace estimation for multivariate observations of high
  dimension: the deterministic signals case.
\newblock {\em IEEE Trans. Inform. Theory}, 58(2):1043--1068, 2012.

\bibitem{X}
J.-s. Xie.
\newblock The convergence on spectrum of sample covariance matrices for
  information-plus-noise type data.
\newblock {\em Appl. Math. J. Chinese Univ. Ser. B}, 27(2):181--191, 2012.

\end{thebibliography}
\end{document}